               \def\version{September 4, 2006}                               %

\documentclass[reqno, 11pt]{amsart}

\usepackage{amsmath}
\usepackage{amsfonts}
\usepackage{amssymb}
\usepackage{amsthm}
\usepackage{mathrsfs}
\usepackage{hyperref}

\def\UseSection{
        \numberwithin{equation}{section}
    \theoremstyle{plain}
        \newtheorem{theorem}    {Theorem}[section]
        \DefineTheorems 
}

\def\DefineTheorems{
    
    \newtheorem{lemma}      [theorem] {Lemma}
    
    \newtheorem{prop}       [theorem] {Proposition}
    
    \newtheorem{cor}        [theorem] {Corollary}

    \theoremstyle{definition}
    \newtheorem{defn}       [theorem] {Definition}

    \theoremstyle{definition}

}

\newcommand{\bt}   {\begin{theorem}}
\newcommand{\et}   {\end  {theorem}}
\newcommand{\bl}   {\begin{lemma}}
\newcommand{\el}   {\end  {lemma}}
\newcommand{\bp}   {\begin{prop}}
\newcommand{\ep}   {\end  {prop}}
\newcommand{\bc}   {\begin{cor}}
\newcommand{\ec}   {\end  {cor}}
\newcommand{\bd}   {\begin{defn}}
\newcommand{\ed}   {\end  {defn}}

\newcommand{\ba}   {\begin{array}}
\newcommand{\ea}   {\end  {array}}
\newcommand{\be}   {\begin{enumerate}}
\newcommand{\ee}   {\end  {enumerate}}
\newcommand{\bi}   {\begin{itemize}}
\newcommand{\ei}   {\end  {itemize}}

\def\eq#1\en{\begin{equation}#1\end{equation}}
\def\eqsplit#1\ensplit{
    \begin{equation}\begin{split}#1\end{split}\end{equation}
    }
\def\eqalign#1\enalign{
    \begin{align}#1\end{align}
    }
\def\eqmul#1\enmul{
    \begin{multline}#1\end{multline}
    }
\newcommand{\eqarrstar} {\begin{eqnarray*}}
\newcommand{\enarrstar} {\end{eqnarray*}}
\newcommand{\eqarray}   {\begin{eqnarray}}
\newcommand{\enarray}   {\end{eqnarray}}

\newcommand{\lbeq}[1]  {\label{e:#1}}
\newcommand{\refeq}[1] {\eqref{e:#1}}    

\makeatletter
\newcommand{\labelcounter}[2]{{%
    \stepcounter{#1}
    \protected@write\@auxout{}%
    {\string\newlabel{#2}{{\csname the#1\endcsname}{\thepage}}}%
    {\ref{#2}}
    }}
\makeatother
\newcommand{\sss}   { \scriptscriptstyle }

\newcommand{\Ccal}   {\mathcal{C}}

\newcommand{\Mcal}   {\mathcal{M}}

\newcommand{\Ocal}   {\mathcal{O}}

\newcommand{\spose}[1] {{\hbox to 0pt{#1\hss}} }
\newcommand{\ltapprox} {\mathrel{\spose{\lower 3pt\hbox{$\mathchar"218$}}
 \raise 2.0pt\hbox{$\mathchar"13C$}}}
\newcommand{\gtapprox} {\mathrel{\spose{\lower 3pt\hbox{$\mathchar"218$}}
 \raise 2.0pt\hbox{$\mathchar"13E$}}}

\newcommand{\nin}  {{ \not\in }}

\UseSection   

\newtheorem{bem}[theorem] {Remark}

\newcommand{\CC}{\mathbb{C}}
\newcommand{\RR}{\mathbb{R}}
\newcommand{\R}{\mathbb{R}}
\newcommand{\N}{\mathbb{N}}
\newcommand{\NN}{\mathbb{N}}
\newcommand{\E}{\mathbb{E}}

\renewcommand{\P}{\mathbb{P}}
\renewcommand{\i}{{\rm{i}}}
\newcommand{\sign}{{\rm{sign}}}

\def\1{{\mathchoice {1\mskip-4mu\mathrm l}      
{1\mskip-4mu\mathrm l}
{1\mskip-4.5mu\mathrm l} {1\mskip-5mu\mathrm l}}}
\newcommand{\ssup}[1] {{\scriptscriptstyle{({#1}})}}

\newcommand{\heap}[2]{\genfrac{}{}{0pt}{}{#1}{#2}}

\newenvironment{proofsect}[1]
{\vskip0.1cm\noindent{\bf #1.}\hskip0.5cm}{\vspace{0.3cm}}

\renewenvironment{proof}[1]
{\vskip0.1cm\noindent{\it #1.}\hskip0.5cm}{\vspace{0.3cm}}

\newenvironment{step}[1]{\vskip0.1cm\noindent{\bf Step {#1}: \it}{\it}}

\newcommand{\EE}[1]{\mathbb{E}_{a}\left(#1\right)}
\newcommand{\ind}[1]{1_{#1}}
\newcommand{\pair}[1]{\langle#1\rangle}
\newcommand{\phibar}{\overline{\phi}}

\newcommand{\cofactor}[1]{\operatorname{det}_{#1}}
\newcommand{\expect}  {{\mathbb E}}
\newcommand{\detrm} {{\rm det}}
\newcommand{\Z}{\mathbb{Z}}
\newcommand{\sT}{{\scriptscriptstyle T}}
\newcommand{\sR}{{\scriptscriptstyle R}}

\newcommand{\sS}{{\scriptscriptstyle S}}
\renewcommand{\d}{{\rm d}}

\newcommand{\eps}{\varepsilon}

\newcommand{\supp}{{\operatorname {supp}}}

\usepackage[usenames]{color}

\newcommand{\ch}[1]{{ #1}}
\newcommand{\chdb}[1]{{ #1}}
\newcommand{\chwk}[1]{{ #1}}

\newcommand{\Switchshort}[1]{}    
\newcommand{\Switchlong}[1]{#1}   

\begin{document}

\title[Joint density for local times]
{\large Joint density for the local times\\ of continuous-time Markov chains\Switchlong{:\\
\ch{Extended version}}}

\author[David Brydges, Remco van der Hofstad and Wolfgang K\"onig]{}
\maketitle

\centerline{\sc David Brydges\footnote{Department of Mathematics,
University of British Columbia,
Vancouver, BC V6T 1Z2, Canada.
{\tt db5d@math.ubc.ca}},
Remco van der Hofstad\footnote{Department of Mathematics and Computer Science,
Technical University Eindhoven, Post Box 513, 5600 MB Eindhoven, The
Netherlands. {\tt rhofstad@win.tue.nl}} and
Wolfgang K\"onig\footnote{Mathematisches Institut, Universit\"at Leipzig,
Augustusplatz 10/11, D-04157 Leipzig, Germany.
{\tt koenig@math.uni-leipzig.de}}}

\vskip.5cm

\centerline{\small\version}

\vskip.5cm


\begin{abstract}
    We investigate the local times of a
    continuous-time Markov chain on an arbitrary discrete state space.
    For fixed finite range of the Markov chain,
    we derive an explicit formula for the joint density of all local times
    on the range, at any fixed time. We use standard tools from
    the theory of stochastic processes and finite-dimensional complex calculus.

    We apply this formula in the following directions:
    (1) we derive large deviation upper estimates for the normalized local times
    beyond the exponential scale, (2)
    we derive the upper bound in Varadhan's \chwk{l}emma for any measurable
    functional of the local times, \ch{and} (3) we derive large deviation upper
    bounds for continuous-time simple random walk on large subboxes of
    $\Z^d$ tending to $\Z^d$ as time diverges. \ch{We finally discuss
    the relation of our density formula to the Ray-Knight theorem
    for continuous-time simple random walk on $\Z$, which is analogous to
    the well-known Ray-Knight description of Brownian local times.}
    \Switchlong{In this extended version, we prove that the Ray-Knight
    theorem follows from our density formula.}
    \end{abstract}

\section{Introduction} \lbeq{sec-intro} Let $\Lambda$ be a finite or
countably infinite set and \chdb{let} $A=(A_{x,y})_{x,y\in\Lambda}$ \chdb{be} the generator,
sometimes called the $Q$-matrix, of a  continuous-time
Markov chain $(X_t)_{t\in[0,\infty)}$ on $\Lambda$.
Under the measure $\P_a$, the chain starts at
$X_0=a\in\Lambda$, and by $\E_a$ we denote the corresponding
expectation. The main object of our study are the {\it local times},
{ defined by}
    \eq\lbeq{loctim}
    \ell_{\sT}(x) =\int_{0}^T \1_{\{X_s=x\}}\,\d s,\qquad
    x\in\Lambda,T>0,
    \en
which register the amount of time the chain
spends in $x$ up to time $T$. We have $\langle
\ell_{\sT},V\rangle=\int_0^T V(X_s)\,\d s$ for any bounded function
$V\colon \Lambda\to\R$, where $\langle \cdot,\cdot\rangle$ denotes
the standard inner product on $\R^\Lambda$. We conceive the normalized
local times tuple, $\frac 1T\ell_{\sT}=(\frac 1T\ell_{\sT}(x))_{x\in\Lambda}$,
as a random element of the set $\Mcal_1(\Lambda)$ of probability
measures on $\Lambda$.

The local times tuple $\ell_{\sT}=(\ell_{\sT}(x))_{x\in\Lambda}$, and
in particular its large-$T$ behavior, are of fundamental interest in
many branches and applications of probability theory. We are
particularly interested in the large deviation of $\frac 1T\ell_{\sT}$.
A by now classical result \cite{Ga77, DV75-83} states, for a {\it finite\/}
state space $\Lambda$, a
large deviation principle for $\frac 1T\ell_{\sT}$,
for any starting point $a\in\Lambda$, on
the scale $T$. More precisely, for any closed set $\Gamma\subseteq\Mcal_1(\Lambda)$,
    \eq\lbeq{upperLDP}
    \limsup_{T\to\infty}\frac1 T\log\P_a({\scriptstyle{\frac1 T}}\ell_{\sT}\in \Gamma)
    \leq -\inf_{\mu\in \Gamma}I_A(\mu),
    \en
    and, for any open set $G\subseteq\Mcal_1(\Lambda)$,
    \eq\lbeq{lowerLDP}
    \liminf_{T\to\infty}\frac 1 T\log\P_a({\scriptstyle{\frac1 T}}
    \ell_{\sT}\in G)\geq -\inf_{\mu\in G}I_A(\mu).
    \en
The rate function $I_A$ may be written
    \begin{equation}\label{ratefunction}
    I_A(\mu)=
    -\inf\Big\{\Big\langle Ag,\frac\mu{g}\Big\rangle\,\Big|\, g\colon\Lambda\to(0,\infty)\Big\}.
    \end{equation}
In case that $A$ is a symmetric matrix, $I_A(\mu)=\|(-A)^{\frac 12}\sqrt{\mu}\|_2^2$ is equal to the Dirichlet form of $A$ \ch{applied to $\sqrt{\mu}$}.
The topology used on $\Mcal_1(\Lambda)$ is the weak topology induced by convergence of
integrals against all bounded functions $\Lambda\to\R$, i.e.,
the standard topology of pointwise convergence since $\Lambda$ is assumed finite.
For infinite $\Lambda$, versions of this large deviations principle may be formulated
for the restriction of the chain to some finite subset of $\Lambda$.
A standard way of proving the above principle of large deviations
is via the G\"artner-Ellis theorem; see \cite{DZ98} for more background on large deviation
theory. One of the major corollaries is Varadhan's \chwk{l}emma, which
states that
\begin{equation}\lbeq{Varadhan}
\lim_{T\to\infty}\frac 1 T\log\E_x\big[{\rm e}^{T F(\frac 1T\ell_{\sT})}]=-\inf_{\mu\in\Mcal_1(\Lambda)}\big[I_A(\mu)-F(\mu)\big],
\end{equation}
for any function $F\colon\Mcal_1(\Lambda)\to\R$ that is bounded and
continuous in the above topology. We would like to stress that in many
situations it is the upper bound in \refeq{Varadhan} that is difficult to
prove since $F$ often fails to be upper semicontinuous. (However, often
$F$ turns out to be lower semicontinuous or well approximated by lower
semicontinous functions, so that the proof of the lower bound in \refeq{Varadhan}
is often { simpler}.)

In the present paper, we considerably strengthen the above large deviation
principle and the assertion in \refeq{Varadhan} by presenting an {\it explicit
density of the random variable\/} $\ell_{\sT}$, i.e., a joint density of the
tuple $(\ell_{\sT}(x))_{x\in\Lambda}$, for any fixed $T>0$. We do this for
either a finite state space $\Lambda$ or for the restriction to a finite
subset. This formula opens up several new possibilities, such as
\begin{enumerate}
\item more precise asymptotics for the probabilities in \refeq{upperLDP}
and \refeq{lowerLDP} and for the expectation on the left of \refeq{Varadhan},

\item the validity of  \refeq{Varadhan} for many discontinuous functions $F$,

\item versions of the large deviation principle for rescaled versions
of the local times on  state spaces $\Lambda=\Lambda_{\sT}$ coupled
with $T$ and growing to some infinite set,

\end{enumerate}

Clearly, a closed analytical formula for the density of the local times
is quite interesting in its own right.
Unfortunately, our expression for the local times density is rather involved
and is quite hard to evaluate asymptotically. Actually, not even the nonnegativity of the density
can be easily seen from our formula. Luckily, {\it upper bounds} on the
density  are more easily obtained. We will be able to use these upper bounds
to derive proofs of \refeq{upperLDP} and of the upper
bound in \refeq{Varadhan} for {\it every measurable set\/} $\Gamma$, respectively,
for {\it every measurable function\/} $F$, which is a great improvement.

This paper is organized as follows. In Section \ref{sec-densloctime},
we identify the density of the local times in Theorem
\ref{thm-locdens}, and prove Theorem \ref{thm-locdens}.
In Section \ref{sec-DV}, we use Theorem \ref{thm-locdens} to prove large
deviation upper bounds in Theorem \ref{thm-DV}.
Finally, we close in Section \ref{sec-dis} by discussing our results,
by relating them to the history of the problem \ch{and by discussing
the relation to the Ray-Knight theorem.}

\section{Density of the local times}
\label{sec-densloctime}

In this section, we present our fundamental result, Theorem
\ref{thm-locdens}, which is the basis for everything that follows. By
\eq
R_{\sT}=\supp(\ell_{\sT})=\{X_s\colon s\in[0,T]\}\subseteq\Lambda
\en
we denote the {\it range\/} of the Markov chain. Note that given
$\{R_{\sT}\subseteq R\}$ for some finite set $R\subseteq\Lambda$, the random
tuple $(\ell_{\sT}(x))_{x\in R}$ does not have a density
\ch{with respect to the Lebesgue measure}, since the
event $\{\ell_{\sT}(x)=0\}$ occurs with positive probability for any
$x\in R$, except for the initial site of the chain. However, given
$\{R_{\sT}=R\}$ for some $R\subseteq\Lambda$,
the tuple $(\ell_{\sT}(x))_{x\in R}$ takes values in the simplex
    \eq
    \Mcal^+_{\sT}(R)=\Big\{l\colon R\to(0,\infty)\,
    \big|\,\sum_{x\in R}l(x)=T\Big\},
    \en
which is a convex open subset of the hyperplane in $\R^R$ that is
perpendicular to $\1$. It will turn out that on  $\{R_{\sT}=R\}$,
the tuple $(\ell_{\sT}(x))_{x\in R}$ has
a density with respect to the Lebesgue
measure $\sigma_\sT$ on $\Mcal^+_{\sT}(R)$ defined by
the disintegration of Lebesgue measure into surface measures,
\begin{equation}\label{eq:dsigma}
    \int \d^{R}l \, F (l)
=
    \int_{0}^{\infty} \d T \,\,
    \int_{\Mcal^{+}_{\sT}}\sigma_{\sT} (\d l)\,\, F (l),
\end{equation}
\ch{where $F\colon(0,\infty)^R\to\R$ is bounded and continuous with compact support.}

We need some notation.  Let $R\subseteq \Lambda$ and let
$a,b\in R$.  For a matrix
$M=(M_{x,y})_{x,y\in \Lambda}$
we denote by $\det^{\ssup{R}}_{ab}(M)$ the \chdb{$(b,a)$ cofactor \chwk{of the $R\times R$-submatrix of $M$},
namely the determinant of the
matrix $(1_{x\not =b}M_{x,y}1_{y\not =a}+1_{x=b,y=a})_{x,y\in R}$}.  We write $\det_{ab}$ instead of
$\det^{\ssup{\Lambda}}_{ab}$ when no confusion can arise.  By
$\partial_{l}$ we denote the $\Lambda\times\Lambda$-diagonal matrix
with $(x,x)$-entry $\ch{\partial_{l_x}}$, which is the partial
derivative with respect to $l_x$. Hence,
$\det^{\ssup{R}}_{ab}(M+\partial_{l})$ is a linear differential
operator of order $\ch{|R|-2+\delta_{a,b}}$.

Then our main result reads as follows:

\begin{theorem}[Density of the local times]
\label{thm-locdens}Let $\Lambda$ be a finite or countably
infinite set with at least two elements and let
$A=(A_{x,y})_{x,y\in\Lambda}$ be the
conservative generator
of a continuous-time Markov chain on $\Lambda$.
Fix a finite subset $R$ of $\Lambda$ and sites $a,b\in R$.
Then, for every $T>0$ and for every bounded measurable function
$F\colon \Mcal^+_{\sT}(R)\to \R$,
    \eq \lbeq{density1}
    \E_a\big[F (\ell_{\sT})\1_{\{X_{\sT}=b\}}\1_{\{R_{\sT}=R\}}\big]
    =
    \int_{\Mcal^+_{\sT}(R)} F(l) \rho^{\ssup R}_{ab}(l) \,\sigma_\sT(\d l),
    \en
where, for $l\in \Mcal^+_{\sT}(R)$,
    \eq\lbeq{density2b}
    \rho^{\ssup R}_{ab}(l)
    =
    \cofactor{ab}^{\ssup{R}} \big(-A+\partial_{l}\big)
    \int_{[0,2\pi]^R}
    {\rm e}^{\sum_{x,y\in R} A_{x,y}\sqrt{l_x}\sqrt{l_y}{\rm e}^{\i(\theta_{x}-\theta_{y})}}
    \prod_{x\in R}\frac{\d \theta_x}{2\pi}.
    \en
\end{theorem}

Alternative expressions for the density $\rho^{\ssup R}_{ab}$
are found in Proposition~\ref{prop-altern} below.
Note that the density $\rho^{\ssup R}_{ab}$ does not depend on the
values of the generator outside $R$, nor on $T$. The formula for the
density is explicit, but quite involved, in particular as it involves
determinants of large matrices, additional multiple integrals, and
various partial derivatives. For example, it is not clear from
\chdb{\refeq{density2b}} that $\rho^{\ssup R}_{ab}$ is
non-negative. Nevertheless, the
formula allows us to prove rather precise and transparent large deviation upper
bounds for the local times \chdb{as we shall see later}.
As we will discuss in more detail in Section \ref{sec-dis},
\chdb{\chwk {T}heorem \ref{thm-locdens} finds its roots in the work of Luttinger
\cite{Lut83} who expressed expectations of functions of the local
times in terms of integrals in which there are ``functions'' of
anticommuting differential forms (Grassman variables).  It is not \chwk{clear from}
his work that the Grassman variables can be removed without
creating intractable expressions. Theorem \ref{thm-locdens}
accomplishes this removal. We also provide a proof that makes no overt
use of Grassman variables; the determinant is their legacy.}


\vskip0.3cm

To prepare for the proof, we need the following two lemmas and some
notation.  We write $\phi = u+\i v$ and $\overline\phi =u-\i v$, where
$u,v\in\RR^{\Lambda}$, and we use $\d^{\Lambda}u
\, \d^{\Lambda}v$ to denote the Lebesgue measure on
$\RR^{\Lambda}\times \RR^{\Lambda}$. Let $\pair{\phi ,\psi}=
\sum_{x\in \Lambda}\phi_{x}\psi_{x}$ be the real inner product on
$\CC^{\Lambda}$.

\begin{lemma}\label{lem:detM} Let $\Lambda$ be a finite set,
and let $M\in\CC^{\Lambda\times\Lambda}$. If $\Re \pair{\phi ,M\overline\phi}>0$ for any $\phi \in\CC^{\Lambda}\setminus\{0\}$,
then
    \begin{equation}\label{eq:detM-1}
    \int_{\RR^{\Lambda}\times
    \RR^{\Lambda}} \d^{\Lambda}u \,\d^{\Lambda}v \,{\rm e}^{- \pair{\phi ,M\phibar}}
    =
        \frac{\pi^{|\Lambda|}} {\detrm( M)}.
    \end{equation}
\end{lemma}

\begin{bem}\label{notation} By introducing polar coordinates $(l,\theta)\in
[0,\infty)^{\Lambda}\times [0,2\pi]^{\Lambda}$ via
    \eq
    \lbeq{polar}
    { \phi_{x} = \sqrt{l_{x}}{\rm e}^{\i\theta_{x}}},\qquad x\in\Lambda,
    \en
we can transform
    \begin{equation}\label{eq:polar1}
    \d^{\Lambda}u\,\d^{\Lambda}v =
    \pi^{|\Lambda|} \prod_{x\in \Lambda} \Big(\d l_{x} \,
    \frac{\d\theta_x}{2\pi}\Big)=2^{-|\Lambda|}\d^\Lambda l\,\d^\Lambda \theta
    \end{equation}
and can rewrite \eqref{eq:detM-1} in the form
    \begin{equation}\label{eq:detM-2}
    \int_{[0,\infty)^{\Lambda}\times [0,2\pi]^{\Lambda}} \d^\Lambda l\,\frac{\d^\Lambda
    \theta}{(2\pi)^{|\Lambda|}}\,
    {\rm e}^{- \pair{\phi ,M\overline\phi}}=\frac 1{\detrm(M)}.
    \end{equation}
\hfill$\Diamond$
\end{bem}

\begin{proofsect}{Proof of Lemma~\ref{lem:detM}}\label{proof-lemma-detM}
We define the complex inner product $ (\phi ,\psi) = \pair{\phi
,\overline{\psi}}$. Any unitary matrix $U\in\CC^{\Lambda\times\Lambda}$ defines a
complex linear transformation on $\CC^{\Lambda}$ by $\phi ' = U\phi$.
By writing $\phi = u+\i v$ and $\phi'=u'+\i v'$ we obtain a real linear
transformation $\widetilde{U}\colon (u,v)\mapsto (u',v')$ on
$\RR^{\Lambda}\oplus \RR^{\Lambda}$. The map $\widetilde{U}$ is orthogonal, because
    \[
    \pair{u',u'}+\pair{v',v'}
    =
    (\phi' ,\phi')
    =
    (U\phi ,U\phi)
    =
    (\phi ,\phi)
    =
    \pair{u,u}+\pair{v,v}.
    \]

Let $M^{\ast}$ be the adjoint to $M$ so that $(\phi
,M\psi)=(M^{\ast}\phi,\psi)$. First we consider the case where $M=M^{\ast}$. The hypothesis
{ $\Re \pair{\phi ,M\overline\phi}>0$} can be rewritten as $(\phi ,M\phi)>0$,
so that $M$ has throughout positive eigenvalues $\lambda_{x}$, $x\in\Lambda$.
Since $M$ is self-adjoint there exists a unitary transformation $U$
such that $U^{\ast}MU = D$, where $D$ is diagonal with diagonal entries
$D_{x,x}=\lambda_{x}>0$. Thus, by the change of variables
$(u',v')=\widetilde{U} (u,v)$,
    \[
    \int \d^{\Lambda}u \,\d^{\Lambda}v \,{\rm e}^{- (\phi ,M\phi)}
    =
    \int \d^{\Lambda}u \,\d^{\Lambda}v \,{\rm e}^{- (\phi ,D\phi)}.
    \]
The integral on the right hand side factors into a product of integrals
    \[
    \prod_{x\in\Lambda}
    \int_{\RR} \d u\,\int_{\RR} \d v \,{\rm e}^{- \lambda_{x} u^{2}- \lambda_{x} v^{2}}
    =
    \prod_{x\in\Lambda}\frac\pi{\lambda_{x}}
    =
    \frac{\pi^{|\Lambda|}}{ \detrm (M)}.
    \]
The lemma is proved for the case $M=M^{\ast}$.

Now we turn to the case where $M^{\ast}\not =M$. Let
    \[
    S = \frac{1}{2}\big(M+M^{\ast}\big) \qquad\mbox{and}\qquad
    A = \frac{1}{2\i}\big(M-M^{\ast}\big).
    \]
Thus, $S$ and $A$ are self-adjoint and $M=S+\i A$. Also, $(\phi ,S\phi)= \Re
\pair{\phi ,M\phibar}$ which is positive by the hypothesis. Therefore
the eigenvalues of $S$ are strictly positive.

For ${ \mu} \in \CC$ we define $ M (\mu) = S+\mu A$.  For
$\mu$ real, the matrix $M (\mu)$ is self-adjoint.
Observe that $M (\mu)$ has throughout strictly positive eigenvalues when
$\mu =0$. Hence, the real part of the characteristic polynomial of $M(\mu)$ is non-zero on $(-\infty,0]$, and therefore bounded away from zero on $(-\infty,0]$, for $\mu=0$. By continuity of the real part of this polynomial in $\mu$, the latter property persists to all $\mu $ in a suitable open interval $I\subset\RR$
containing the origin. Therefore, $ M (\mu)$ has throughout
strictly positive eigenvalues for all $\mu\in I$.
Thus we have $(\phi ,M (\mu)\phi)>0$ for all nonzero $\phi$ and all $\mu\in I$.

Now we apply the preceding with $M=M (\mu)$, and obtain, for
$\mu \in I$,
    \begin{equation}\label{eq:detM-3}
    \detrm (M (\mu))
    \int \d^{\Lambda}u \,\d^{\Lambda}v \,{\rm e}^{- (\phi ,M (\mu)\phi)}
    =
    \pi^{|\Lambda|}.
    \end{equation}
Both sides of this equation are analytic in $\mu$ for $\Re \mu
\in I$ because $\detrm (M (\mu))$ is a polynomial in $\mu$, and
the integral of the analytic function $\exp (- (\phi ,M
(\mu)\phi))$ is analytic by Morera's theorem and the Fubini
theorem, as well as the remark that
    \[
    \big|{\rm e}^{- (\phi ,M (\mu)\phi)}\big|
    =
    \big|{\rm e}^{- (\phi ,S\phi) - \mu (\phi ,A\phi)}\big|
    =
    {\rm e}^{- (\phi ,S\phi) - \Re\mu (\phi ,A\phi)}
    =
    {\rm e}^{- (\phi ,M (\Re \mu)\phi)}.
    \]
By analytic continuation (\ref{eq:detM-3}) holds for $\Re \mu \in
I$ and in particular for $\mu =\i$. At $\mu =\i$, $M
(\mu)=M$. \qed
\end{proofsect}
\vskip0.3cm

\begin{lemma}\label{lem:conjugate}
Let $\Lambda$ be a finite set, let $M\in\CC^{\Lambda\times\Lambda}$, and $v=(v_x)_{x\in\Lambda}\in \CC^{\Lambda}$. Then, for any continuously differentiable function $g\colon\CC^{\Lambda}\to\R$,
\eq\label{conjugate-1}
    \cofactor{ab} (M+\partial_{l})\big(
    {\rm e}^{\langle v,\cdot\rangle} g\big)(l)
=
    {\rm e}^{\langle v,l\rangle}
    \cofactor{ab}(M+V+\partial_{l})g(l),\qquad l\in\RR^\Lambda,
\en
where $V=(\delta_{xy}v_x)_{x,y\in\Lambda}$ denotes the diagonal matrix with diagonal entries $v_x$.
\end{lemma}

\begin{proofsect}{Proof}
By a cofactor expansion, one sees that, for any diagonal matrix $W$, $\cofactor{ab} (M+W) =
\sum_{X\subseteq \Lambda\setminus \{a,b \}}c_{X} \prod_{x\in X}W_{x,x}$ for suitable coefficients $c_X$ depending only on the entries of $M$. Analogously, $\cofactor{ab} (M+\partial_{l}) =
\sum_{X\subseteq \Lambda\setminus \{a,b \}}c_{X} \partial_{l}^{X}$, where we used the notation $\partial_{l}^{X}= \prod_{x\in
X} \partial_{l_{x}}$. Therefore,
\begin{align*}
    {\rm e}^{-\langle v,l\rangle}\cofactor{ab} (M+\partial_{l})\big({\rm e}^{\langle v,\cdot\rangle} g\big)(l)
&=
    \sum_{X\subseteq \Lambda\setminus \{a,b \}} c_{X}
    {\rm e}^{-\langle v,l\rangle}
    \partial_{l}^{X}
    \big({\rm e}^{\langle v,\cdot\rangle}g\big)(l)\\
&=
    \sum_{X\subseteq \Lambda\setminus \{a,b \}} c_{X}
    \prod_{x\in\Lambda}\big(v_x +\partial_{l_x}\big){\chdb{g (l)}}
=
    \cofactor{ab}(M+V+\partial_{l})g(l).
\end{align*}
\qed
\end{proofsect}

\begin{proofsect}{Proof of Theorem~\ref{thm-locdens}} We have divided the proof into six steps. In the first
five steps we assume that $\Lambda$ is a finite set, and we put
$R=\Lambda$. { Recall the notation in Remark~\ref{notation}, which will be used
throughout this proof.} We abbreviate $D_l =
\cofactor{ab}(-A-\partial_{l})$.

\begin{step}{1} For any $v\in\CC^{\Lambda}$ with $\Re v\in
(-\infty,0)^\Lambda$, for $F(l) = {\rm e}^{\langle v,l\rangle}$,
    \eq
    \lbeq{baseq}
    \int_0^{\infty} \expect_a \big[
   F(\ell_{\sT}) \1_{\{X_{\sT}=b\}}\big]\,\d T
    = \int_{[0,\infty)^\Lambda\times [0,2\pi]^\Lambda}  (D_lF) (l) {\rm e}^{\langle \phi,A\phibar\rangle} \,\,\prod_{x\in \Lambda} \Big(\d l_{x} \,
    \frac{\d\theta_x}{2\pi}\Big).
    \en
\end{step}

\begin{proof}{Proof}
Recall that $\langle v,
\ell_{\sT}\rangle=\int_0^T v(X_s)\,\d s$ to obtain
    \eq\lbeq{startF1}
    \begin{aligned}
    \int_0^{\infty} \expect_a \big[
   F(\ell_{\sT}) \1_{\{X_{\sT}=b\}}\big]\,\d T
    &=\int_0^{\infty}\expect_a \Big[
    {\rm e}^{\int_0^T v(X_s)\,\d s}\1_{\{X_{\sT}=b\}}\Big]\,\d T=\int_0^{\infty}  \big({\rm e}^{T(A+V)}\big)_{a,b}\,\d T\\
    &=(-A-V)^{-1}_{a,b},
    \end{aligned}
    \en
where $V$ is the diagonal matrix with $(x,x)$-entry $v_x$, and $M_{x,y}$ denotes
the $(x,y)$-entry of a matrix $M$. \ch{In order to see the last identity in
\refeq{startF1}, we note that
    \eq
    \int_0^{\infty}  \big({\rm e}^{T(A+V)}\big)_{a,b}\,\d T
    =\Big(\int_0^{\infty}  {\rm e}^{T(A+V)}\,\d T\Big)_{a,b},
    \en
and that
    \eq
    (A+V)\int_0^{\infty}{\rm e}^{T(A+V)} \,\d T
    = \int_0^{\infty} \frac{\chwk{\d}}{\chwk{\d}T} {\rm e}^{T(A+V)} \,\d T
    =-I.
    \en
}

By Cramer's rule followed by
(\ref{eq:detM-2}),
  \eq\label{eq:detM-2a}
    \begin{aligned}
    (-A-V)^{-1}_{a,b}
    =\frac{\cofactor{ab}(-A-V)}{ \det (-A-V)}
    &=\int  \cofactor{ab}(-A-V) {\rm e}^{\langle \phi,(A+V)\phibar\rangle} \,\,\prod_{x\in \Lambda}
    \Big(\d l_{x} \, \frac{\d\theta_x}{2\pi}\Big)\\
    &=\int  \cofactor{ab}(-A-V) {\rm e}^{\langle v,l\rangle} {\rm e}^{\langle \phi,A\phibar\rangle} \,\,\prod_{x\in \Lambda}
    \Big(\d l_{x} \, \frac{\d\theta_x}{2\pi}\Big).
    \end{aligned}
    \en
We use Lemma~\ref{lem:conjugate} with $g=1$ and $M=A$ to obtain that
    \eq
    \begin{aligned}
    \cofactor{ab}(-A-V) {\rm e}^{\langle v,l\rangle}
    &= (-1)^{|\Lambda|-1}{\rm e}^{\langle v,l\rangle}\cofactor{ab}(A+V)
    =(-1)^{|\Lambda|-1}\cofactor{ab}(A+\partial_l) {\rm e}^{\langle v,l\rangle}\\
    &=\cofactor{ab}(-A-\partial_l) {\rm e}^{\langle v,l\rangle}=(D_l F)(l),
    \end{aligned}
    \en
    where we recall that $D_l =
\cofactor{ab}(-A-\partial_{l})$. Substituting this in \eqref{eq:detM-2a} and combining this with  \refeq{startF1}, we
conclude that \refeq{baseq} holds.
\qed\end{proof}

\begin{step}{2} The formula \refeq{baseq} \ch{is also valid for} functions $F$ of the form
    \eq
    \lbeq{F-def}
    F(l)=
    \prod_{x\in \Lambda}\big( {\rm e}^{v_{x}l_{x}}f_{x}(l_x)\big),\qquad f_x\in\Ccal^2\big((0,\infty)\big),\supp(f_x)\subseteq (0,\infty)\mbox{ compact}, \Re v_x<0.
    \en
    \end{step}

\begin{proof}{Proof}  Note that \refeq{baseq} is linear in $F$ and so if we
know it for exponentials, { then} we obtain it for linear combinations of
exponentials. In more detail, consider the Fourier representation
$f_{x}(l_x) = \int_\RR \widehat{f}_{x} (w_{x}) {\rm e}^{\i w_{x}l_{x}}\,
 \d w_{x}$. Apply \refeq{baseq} for $v$ replaced by $v+\i w$ with $w\in\RR^\Lambda$ to obtain
 $$
 \int_0^{\infty} \expect_a \big[
   {\rm e}^{\langle v,\ell_\sT\rangle}{\rm e}^{\i\langle w,\ell_\sT\rangle} \1_{\{X_{\sT}=b\}}\big]\,\d T
    = \int_{[0,\infty)^\Lambda\times [0,2\pi]^\Lambda}  \big(D_l{\rm e}^{\langle v+\i w,l\rangle}\Big) {\rm e}^{\langle \phi,A\phibar\rangle} \,\,\prod_{x\in \Lambda} \Big(\d l_{x} \,
    \frac{\d\theta_x}{2\pi}\Big).
    $$
Now multiply both sides with $\prod_{x\in\Lambda}\widehat{f}_{x} (w_{x})$ and integrate over $\RR^\Lambda$ with respect to $\d^{\Lambda}w$.  Then we apply Fubini's
theorem to move the $\d^{\Lambda}w$ integration inside. From the representation
 \begin{align*}
    \widehat{f}_{x} (w_x)
=
    \frac{1}{2\pi }\int_\RR f_{x} (l_x) {\rm e}^{-\i w_x l_x}\,\d l_x
\end{align*}
we see that $\widehat{f}_{x}$ is continuous by the dominated
convergence theorem.  Furthermore, $\widehat{f}_{x}$ satisfies the bound
\begin{align*}
    |\widehat{f}_{x} (w_x) |=
    \Big|\frac{1}{2\pi (\i w_x)^{2}}\int f_{x} (l) \frac{\d^{2}}{\d w_x^{2}}{\rm e}^{-\i w_x l}\,\d l\Big|=
    \frac{1}{2\pi }\frac{1}{w_x^{2}}\Big|\int f_{x}''(l) {\rm e}^{-\i w_x l}\,\d l\Big|\leq \frac{1}{2\pi }\frac{1}{w_x^{2}}\int |f_x'' (l) | \,\d l.
\end{align*}
Hence, all functions $w_x\mapsto \widehat{f}_{x} (w_x)$ are absolutely integrable, and the
exponentials with $\Re v_{x}<0$ make the integration over $l_x$
convergent for any $x\in R$.
\qed\end{proof}

In the following we abbreviate
$D_l^{\ast}=\cofactor{ab}(-A+\partial_{l})$.

\begin{step}{3} For $F$ as in \refeq{F-def},
   \eq\lbeq{baseq2}
   \int_0^{\infty} \expect_a \big[
    F(\ell_{\sT})\1_{\{X_{\sT}=b\}}\big]\,\d T
    =\int F(l) \ch{D^{\ast}_{l}}{\rm e}^{\langle \phi,A\phibar\rangle} \,\,\prod_{x\in \Lambda}\Big( \d l_{x} \, \frac{\d\theta_x}{2\pi}\Big).
    \en
 \end{step}
    \begin{proof}{Proof} Comparing \refeq{baseq} with this formula we see that it is
enough to prove that the integration by parts formula
    \eq
    \lbeq{claimPI}
    \int (\ch{D_l}F)(l)\,
    {\rm e}^{\langle \phi,A\phibar\rangle} \,\,\d^{\Lambda}l
    =\int\,  F(l)\,\big(\ch{D^{\ast}_l}{\rm e}^{\langle \phi,A\phibar\rangle}\big)
     \,\,\d^{\Lambda}l
    \en
holds for any  $\theta\in[0,2\pi]^\Lambda$. Since $\ch{D_l}=\cofactor{ab}(-A-\partial_{l})$
is a linear differential operator which is first order in each partial
derivative, it suffices to consider one integral at a time and perform the
integration by parts as follows: for any $x\in\Lambda$ and any fixed $(l_y)_{y\in\Lambda\setminus\{x\}}$,
    \eq\lbeq{f1claim}
    \int_0^{\infty} \big(-\partial_{l_x}F(l)\big){\rm e}^{\langle \phi,A\phibar\rangle}\,\d l_x
    =\int_0^{\infty} F(l) \partial_{l_x} {\rm e}^{\langle \phi,A\phibar\rangle}\,\d l_x,\qquad x\in\Lambda.
    \en
There are no boundary contributions because the map $l_x\mapsto F(l)$  has a compact support in $(0,\infty)$. This proves \refeq{baseq2}.
\qed\end{proof}

\begin{step}{4}For any $v\in\CC^\Lambda$,
    \eq\lbeq{baseq2b}
   \int_0^{\infty} \expect_a \big[
    {\rm e}^{\langle v,\ell_{\sT}\rangle}
    \1_{\{X_{\sT}=b\}}\1_{\{R_{\sT}=\Lambda\}}
    \big]\,\d T
    =\int {\rm e}^{\langle v,l\rangle}
    \ch{D^{\ast}_l}{\rm e}^{\langle \phi,A\phibar\rangle}
    \,\,\prod_{x\in \Lambda} \Big(\d l_{x} \, \frac{\d\theta_x}{2\pi}\Big).
    \en
    \end{step}
    \begin{proof}{Proof}
    Let $(f_{n})_{n\in\NN}$ be a uniformly bounded sequence of smooth functions
with compact support in $(0,\infty)$ such that $f_{n} (t) \rightarrow
\1_{(0,\infty)}(t)$ for any $t$. Choose $F(l)=F_n (l)=\prod_{x\in\Lambda} ({\rm
e}^{v_{x}l_{x}}f_{n} (l_{x}))$ in \refeq{baseq2} and take the limit as
$n\rightarrow \infty$, interchanging the limit with the integrals
using the dominated convergence theorem. Observe that $\lim_{n\to\infty}F_n(\ell_\sT)={\rm e}^{\langle v,\ell_\sT\rangle}\prod_{x\in\Lambda} \1_{(0,\infty)}(\ell_\sT(x))={\rm e}^{\langle v,\ell_\sT\rangle}\1_{\{R_{\sT}=\Lambda\}}$ almost surely. Furthermore, $\lim_{n\to\infty}F_n(l)={\rm e}^{\langle v,l\rangle}$ almost everywhere with respect
to the measure $\prod_{x\in \Lambda} (\d l_{x} \,
\frac{\d\theta_x}{2\pi})$. Thus we obtain \refeq{baseq2b} in the limit of
\refeq{baseq2}.
\qed\end{proof}

\begin{step}{5} For all $v\in \CC^{\Lambda}$,
    \eq
    \lbeq{baseq3}
    \expect_a \big[
    {\rm e}^{\langle v,\ell_{\sT}\rangle}
    \1_{\{X_{\sT}=b\}}\1_{\{R_{\sT}=\Lambda\}}\big]
    =\int_{\Mcal^+_{\sT}(\Lambda)}
    {\rm e}^{\langle v,l\rangle}
    \rho_{ab}^{\ssup\Lambda}(l)\,\d^\Lambda l,\qquad T>0, a,b\in \Lambda,
    \en
where $\rho_{ab}^{\ssup\Lambda}(l)$ is given by
\refeq{density2b}.
\end{step}

    \begin{proof}{Proof} Recall that $  \sum_{x\in\Lambda} \ell_{\sT}(x)=T$ almost surely and that $\sum_{x\in\Lambda} l_{x}=T$ for $l\in\Mcal^+_{\sT}(\Lambda)$. Hence, without loss of generality, we can assume that $\Re v \in
(-\infty ,0)^{\Lambda}$, since adding a constant $C\in \RR$ to all the $v_x$ results in adding a factor of ${\rm e}^{CT}$ on both sides. In \refeq{baseq2b} we replace $v_{x}$ by $v_{x}-\lambda$ with $\lambda >0$. Then
\refeq{baseq2b} becomes
   \eq\lbeq{baseq2c}
   \begin{aligned}
   \int_0^{\infty} {\rm e}^{-\lambda T}\expect_a \big[
   {\rm e}^{\langle v,\ell_{\sT}\rangle}
    \1_{\{X_{\sT}=b\}}\ch{\1_{\{R_{\sT}=\Lambda\}}}\big]\,\d T
    &=\int
    {\rm e}^{\langle v,{ l}\rangle}{\rm e}^{-\lambda \sum_x l_{x}}
    D_l^{\ast}{\rm e}^{\langle \phi,A\phibar\rangle}\,\,\prod_{x\in \Lambda} \Big(\d l_{x} \, \frac{\d\theta_x}{2\pi}\Big)\\
    &=\int_{(0,\infty)^\Lambda}
    {\rm e}^{\langle v,{ l}\rangle}{\rm e}^{-\lambda \sum_x l_{x}}
    \rho_{ab}^{\ssup\Lambda}(l)\,\,\d^{\Lambda} l\\
    &=\int_0^{\infty}{\rm e}^{-\lambda T}
    \Big[\int_{\Mcal^+_{\sT}(\Lambda)}
    {\rm e}^{\langle v,{ l}\rangle}
    \rho_{ab}^{\ssup\Lambda}(l)\,\sigma_{\sT}(\d l)\Big]\,\d T,
    \end{aligned}
    \en
    where
\begin{equation}\label{rhofirst}
\rho_{ab}^{\ssup\Lambda}(l)=\int_{[0,2\pi]^\Lambda}D_l^{\ast}{\rm e}^{\langle \phi,A\overline\phi\rangle}\,\prod_{x\in \Lambda}\frac{\d\theta_x}{2\pi},\qquad l\in (0,\infty)^\Lambda.
\end{equation}
In the second equation, we have interchanged the integrations over $l$ and
$\theta$ and have rewritten the $\theta$ integral using \refeq{density2b}.
In the third equation in \refeq{baseq2c}, we have introduced the
variable $T= \sum_x l_{x}$ and
used (\ref{eq:dsigma}).

\chdb{Hence we have proved that the Laplace transforms with respect to
$T$ of the two sides of \refeq{baseq3} coincide. \chwk{As a consequence, \refeq{baseq3} holds for almost every $T>0$.} Furthermore, \refeq{baseq3}
even holds for all $T>0$, since both sides are continuous. Indeed, for
small $h$ we have $\langle v, \ell_{T+h}\rangle =\langle v,
\ell_T\rangle$, $X_{T+h}=X_T$ and $R_{T+h}=R_T$ with high probability,
which easily implies the continuity of the left hand side of
\refeq{baseq3}. We \chwk{see} that the right hand side is continuous for $T> 0$
by using the change of variable $t=T^{-1}l$ and \refeq{density2b} to
rewrite the right hand side as an integral of a continuous function of
$T,t$ on the standard simplex $\Mcal^+_{\sT=1}(\Lambda)$.}
\qed\end{proof}

Now we { complete} the proof of the theorem:

\begin{step}{6} The formula \refeq{density1} holds for any finite or countably
infinite state space $\Lambda$ and any  finite subset $R$ of $\Lambda$.
\end{step}

    \begin{proof}{Proof}
It is enough to prove \refeq{density1} for the case $F (l)= {\rm
e}^{\langle v,l\rangle}$ with $\Re(v)\in(-\infty,0)^R$ because the distribution of
$(\ell_{\sT} (x))_{x\in R}$ on the event $\{R_{\sT}=R\}$ is determined by
its characteristic function.

Consider the Markov chain on $R$ with conservative generator
$A^{\ssup{R}}=(A^{\ssup{R}}_{x,y})_{x,y\in R}$ given by
    \eq\lbeq{ARdef}
    A^{\ssup R}_{x,y}=\begin{cases}A_{x,y}&\mbox{if } x\neq y,\\
    -\sum_{y\in R\setminus\{ x\}}A_{x,y}&\mbox{if }x=y,
    \end{cases}
    \en
and let $V^{\ssup{R}}$ be the diagonal $R\times R$ matrix with $V^{\ssup{R}}_{x,x}= \sum_{y\in
\Lambda\setminus R}A_{x,y}$. Then
    \eq\lbeq{AR-A}
    { A^{\ssup{R}}_{x,y}
    =
    A_{x,y}+V^{\ssup{R}}_{x,y} \qquad\forall x,y\in R.}
    \en
When started in $R$, the Markov chain { with generator
$A^{\ssup{R}}$} coincides with the original
one as long as no step to a site outside $R$ is attempted. { Step
decisions outside $R$ are suppressed.}
The distribution of this chain is absolutely continuous
with respect to the original one. More precisely,
    \eq
    \lbeq{baseq4}
    \expect_a \Big[
    F(\ell_{\sT})\1_{\{X_{\sT}=b\}}\1_{\{R_{\sT}=R\}}\Big]
    = \expect^{\ssup R}_{a}
    \Big[F(\ell_{\sT})
    {\rm e}^{-\sum_{x\in R}\ell_{\sT}(x)V^{\ssup{R}}_{x,x}}
    \1_{\{X_{\sT}=b\}}\1_{\{R_{\sT}=R\}}\Big],\qquad T>0, a,b\in R,
    \en
where $\expect^{\ssup R}_{a}$ is the expectation with respect to the
Markov chain on $R$ with generator $A^{\ssup R}$. Applying \refeq{baseq3} for this
chain with ${\rm e}^{\langle v,l\rangle}$ replaced by
    \eq
    F_{\sR}(l)=F (l)
    {\rm e}^{-\sum_{x\in R}{ l_x}V^{\ssup{R}}_{x,x}}
    \en
and with $\Lambda$ replaced by $R$, we obtain, writing $\partial^{\ssup
R}_l$ for the restriction of $\partial_{l}$ to $R\times R$,
   \eq\lbeq{densR}
   \begin{aligned}
   \expect_a \Big[
    F(\ell_{\sT})\1_{\{X_{\sT}=b\}}\1_{\{R_{\sT}=R\}}\Big]
    &=\expect^{\ssup R}_{a}
    \Big[F_{\sR}(\ell_{\sT})\1_{\{X_{\sT}=b\}}\1_{\{R_{\sT}=R\}}\Big]\\
    &= \int_{\Mcal^+_{\sT}(R)}F_{\sR}(l) \rho_{ab}^{\ssup{R}}(l)\,\sigma_\sT(\d l)
    = \int_{\Mcal^+_{\sT}(R)}F(l)
    \widetilde\rho_{ab}^{\ssup{\Lambda ,R}}(l)\,\sigma_\sT(\d l)
    \end{aligned}
    \en
where
    \eq\lbeq{rhotildeR}
    \widetilde\rho_{ab}^{\ssup{\Lambda, R}}(l)
    =
    {\rm e}^{-\sum_{x\in R}{ l_x}V^{\ssup{R}}_{x,x}}
    \cofactor{ab}(-A^{\ssup{R}}+\partial^{\ssup{R}}_l)
    \int_{[0,2\pi]^R}
    {\rm e}^{\sum_{x,y \in R} \phi_{x}A^{\ssup{R}}_{x,y}\phibar_{y}}
    \prod_{x\in R} \frac{\d\theta_x}{2\pi}.
    \en
By Lemma~\ref{lem:conjugate} followed by \refeq{AR-A},
    \eq\lbeq{rhotildeR-2}
    \begin{aligned}
    \widetilde\rho_{ab}^{\ssup{\Lambda, R}}(l)
    &=
    \cofactor{ab}(-A^{\ssup{R}}\ch{-}V^{\ssup{R}}+\partial^{\ssup{R}}_l)
    \Big[{\rm e}^{-\sum_{x\in R}l_x V^{\ssup{R}}_{x,x}}
    \int_{[0,2\pi]^R}
    {\rm e}^{\sum_{x,y \in R} \phi_{x}A^{\ssup{R}}_{x,y}\phibar_{y}}
    \prod_{x\in R} \frac{\d\theta_x}{2\pi}\Big]\\
    &=
    \cofactor{ab}(-A^{\ssup{R}}\ch{-}V^{\ssup{R}}+\partial^{\ssup{R}}_l)
    \int_{[0,2\pi]^R}
    {\rm e}^{\sum_{x,y \in R} \phi_{x}A_{x,y}\phibar_{y}}
    \prod_{x\in R} \frac{\d\theta_x}{2\pi}\\
    &=
    \cofactor{ab}^{\ssup{R}}(-A+\partial_{l})
    \int_{[0,2\pi]^R}
    {\rm e}^{\sum_{x,y \in R} \phi_{x}A_{x,y}\phibar_{y}}
    \prod_{x\in R} \frac{\d\theta_x}{2\pi}.
    \end{aligned}
 \en
From the definition \refeq{density2b}, { and using \refeq{polar},}
we recognise the last line as
$\rho^{\ssup R}_{ab}(l)$. Therefore, by combining \refeq{rhotildeR-2}
and \refeq{densR} we have proved \refeq{density1} in the theorem.
\qed
\end{proof}
\end{proofsect}

Now we collect some \ch{alternative} expressions for the density $\rho^{\ssup R}_{ab}$:

\begin{prop}\label{prop-altern}
Let the assumptions of Theorem~\ref{thm-locdens} be satisfied. Let $B=([1-\delta_{x,y}]A_{x,y})_{x,y\in\Lambda}$ be the off-diagonal part of $A$. Then, for any finite subset $R$ of $\Lambda$ and for any sites $a,b\in R$, and for any $l\in \Mcal^+_{\sT}(R)$, the following holds:

\begin{enumerate}
\item[(i)]
\begin{equation}\lbeq{density2}
    \rho^{\ssup R}_{ab}(l)
    = {\rm e}^{\sum_{x\in R} l_x A_{x,x}}
    \cofactor{ab}^{\ssup{R}} \big(-B+\partial_{l}\big)\int_{[0,2\pi]^R}
    {\rm e}^{\sum_{x,y\in R} B_{x,y}\sqrt{l_x}\sqrt{l_y}{\rm e}^{\i(\theta_{x}-\theta_{y})}}
    \prod_{x\in R}\frac{\d \theta_x}{2\pi}.
    \end{equation}

\item[(ii)] For any $r\in(0,\infty)^R$,
\begin{equation}\lbeq{density3}
\rho^{\ssup R}_{ab}(l)= {\rm e}^{\sum_{x\in R} l_x A_{x,x}}
    \cofactor{ab}^{\ssup{R}} \big(-B+\partial_{l}\big)
    \int_{[0,2\pi]^R}
    {\rm e}^{\sum_{x,y\in R} r_{x}B_{x,y}r_{y}^{-1}\,\sqrt{l_x}\sqrt{l_y}{\rm e}^{\i(\theta_{x}-\theta_{y})}}
    \prod_{x\in R}\frac{\d \theta_x}{2\pi}.
    \end{equation}

    \item[(iii)]
\begin{equation}\lbeq{density4}
\rho^{\ssup R}_{ab}(l)
    =
    \int_{[0,2\pi]^R}\cofactor{ab}^{\ssup{R}} \big(-B+V_{\theta,l}\big){\rm e}^{\sum_{x,y\in R} A_{x,y}\sqrt{l_x}\sqrt{l_y}{\rm e}^{\i(\theta_{x}-\theta_{y})}}
    \prod_{x\in R}\frac{\d \theta_x}{2\pi},
    \end{equation}
    where $V_{\theta,l}=(\delta_{x,y}v_{\theta,l}(x))_{x\in R}$ is the diagonal matrix with entries
    \begin{equation}
    v_{\theta,l}(x)=\sum_{z\in R} B_{x,z}\sqrt{\frac{l_z}{l_x}}{\rm e}^{\i(\theta_{x}-\theta_{z})},\qquad x\in R.
    \end{equation}
\end{enumerate}
\end{prop}

The formula in \refeq{density3} will be helpful later when we
derive upper bounds on $\rho^{\ssup R}_{ab}(l)$ in the case
that $A$ is not symmetric. \ch{The remainder of the paper does not
rely on the formula in \refeq{density4}. However, we find
\refeq{density4} of independent interest, since the
integral in \refeq{density4} does not involve any derivative.}

\begin{proofsect}{Proof} Formula \refeq{density2} follows from
\refeq{density2b} by using Lemma~\ref{lem:conjugate}.

We now prove \refeq{density3}. Fix $r\in(0,\infty)^R$ and observe that, for any $l\in(0,\infty)^R$,
\begin{equation}\label{circletrick}
\int_{[0,2\pi]^R}
    {\rm e}^{\sum_{x,y\in R} B_{x,y}\sqrt{l_x}\sqrt{l_y}{\rm e}^{\i(\theta_{x}-\theta_{y})}}
    \prod_{x\in R}\frac{\d \theta_x}{2\pi}
  =\int_{[0,2\pi]^R}
    {\rm e}^{
    \sum_{x,y\in R}
    r_{x}B_{x,y}r_{y}^{-1}\,
    \sqrt{l_x}\sqrt{l_y}{\rm e}^{\i(\theta_{x}-\theta_{y})}
    }
    \prod_{x\in R}\frac{\d \theta_x}{2\pi}.
    \end{equation}
Indeed, substituting ${\rm e}^{\i\theta_{x}}=z_{x}$ for $x\in R$, we can rewrite the integrals as integrals over circles in the complex plane. The integrand is analytic in $z_x\in\CC\setminus\{0\}$. Hence, the integral is independent of the curve (as long as it is closed and winds around zero precisely once), and it is equal to the integral along the centred circle with radius $r_x$ instead of radius one. Re-substituting $r_x {\rm e}^{\i\theta_{x}}= z_{x}$, we arrive at \eqref{circletrick}.
Comparing to \refeq{density2}, we see that we have derived \refeq{density3}.

Finally, we prove \refeq{density4}. We use \refeq{density3} with $r=\sqrt l$ and interchange $\cofactor{ab}^{\ssup{R}} \big(-B+\partial_{l}\big)$ with $\int_{[0,2\pi]^R}$ (this is justified by the analyticity of the integrand in all the $l_x$ with $x\in R$). This gives that
$$
\rho^{\ssup R}_{ab}(l)={\rm e}^{\sum_{x\in R} l_x A_{x,x}}
    \int_{[0,2\pi]^R}\cofactor{ab}^{\ssup{R}} \big(-B+\partial_{l}\big)
    {\rm e}^{\sum_{x,y\in R} l_x B_{x,y}{\rm e}^{\i(\theta_{x}-\theta_{y})}}
    \prod_{x\in R}\frac{\d \theta_x}{2\pi}.
$$
Use Lemma~\ref{lem:conjugate} with $g=1$ to see that
$$
\cofactor{ab}^{\ssup{R}} \big(-B+\partial_{l}\big)
    {\rm e}^{\sum_{x,y\in R} l_x B_{x,y}{\rm e}^{\i(\theta_{x}-\theta_{y})}}
    ={\rm e}^{\sum_{x,y\in R} l_x B_{x,y}{\rm e}^{\i(\theta_{x}-\theta_{y})}}\cofactor{ab}^{\ssup{R}} \big(-B+\widetilde V_\theta\big),
    $$
    where $\widetilde V_\theta=(\delta_{x,y}\widetilde v_\theta(x))_{x\in R}$ is the diagonal matrix with entries $\widetilde v_\theta(x)=\sum_{z\in R} B_{x,z} {\rm e}^{\i(\theta_{x}-\theta_{z})}$.

    Now we use the same transformation as in \eqref{circletrick}: We interpret the integrals over $\theta_x$ as integrals over circles of radius $\sqrt{l_x}$ and replace them by integrals over circles with radius one. By this transformation, $\widetilde V_\theta$ is transformed into $V_{\theta,l}$, and the term ${\rm e}^{\sum_{x,y\in R} l_x B_{x,y}{\rm e}^{\i(\theta_{x}-\theta_{y})}}$ is transformed into ${\rm e}^{\sum_{x,y\in R} \sqrt{l_x} B_{x,y}\sqrt{l_y}{\rm e}^{\i(\theta_{x}-\theta_{y})}}$. Recalling that $B$ is the off-diagonal part of $A$, \refeq{density4} follows.
\qed\end{proofsect}

\section{Large deviation upper bounds for the local times}
\label{sec-DV}

\noindent In this section we use Theorem~\ref{thm-locdens} to derive
sharp upper bounds for the probability in \refeq{upperLDP} and for
the expectation in \refeq{Varadhan} for fixed $T$ and fixed finite
ranges of the local times. The main term in this estimate is given
in terms of the rate function $I_A$. The main value of
our formula, however, comes from the facts that (1) the error term is
controlled on a subexponential scale, (2) the set $\Gamma$ in
\refeq{upperLDP} is just assumed measurable, and (3) the functional
$F$ in \refeq{Varadhan} is just assumed measurable.
Let us stress that this formula is extremely useful, since the functional $F$
is not upper semicontinuous nor bounded in many important
applications.

In Section~\ref{sec-pointwise} we give a pointwise upper bound for the density, in
Section~\ref{sec-uppLDP} we apply it to derive upper bounds for the probability in \refeq{upperLDP} and for the expectation in \refeq{Varadhan}, and in Section~\ref{sec-rescaled}
we consider the same problem for state spaces $\Lambda=\Lambda_{\sT}\subseteq\Z^d$ depending
on $T$ and increasing to $\Z^d$.

\subsection{Pointwise upper bound for the density}\label{sec-pointwise}

\noindent Here is a pointwise upper bound for the density. Recall the rate function $I_A$ introduced in \eqref{ratefunction}.

\begin{prop}[Upper bound for ${\rho_{ab}^{\ssup R}}$]\label{lem-rhobound}
Under the assumptions of Theorem~\ref{thm-locdens}, for any finite
subset $R$ of $\Lambda$, and for any $a,b\in R$, any $T>0$ and any $l\in \Mcal^+_{\sT}(R)$,
\eq\label{rhobound}
    \rho_{ab}^{\ssup R}(l)
\leq {\rm e}^{-TI_A(\frac 1Tl)}
\Big(\prod_{x\in R\setminus \{a,b \}}\sqrt{\frac{T}{l_x}}\Big)\eta_{\sR}^{|R|-1}{\rm e}^{[\chdb{\eta_{\sR}^{-1}}+(4\eta_{\sR}^2 T)^{-1}]\sum_{x,y\in R} \sqrt {l_x}g_y B_{x,y}/(\sqrt{l_y}g_x)},
\en
where $g\in(0,\infty)^R$ is the minimizer in \eqref{ratefunction} \ch{for $\mu=\frac{l}{T}$}
and
    \eq\lbeq{new-CRdef}
    \chdb{\eta_{\sR}}
    =
    \max
    \Big\{
    \max_{x\in R}\sum_{y\in R\setminus \{x\}}\chdb{|B_{x,y}|},\
    \max_{y\in R}\sum_{x\in R\setminus \{y\}}\chdb{|B_{x,y}|},\,1
    \Big\},
    \en
    where $B=([1-\delta_{x,y}]A_{x,y})_{x,y\in\Lambda}$ is the off-diagonal part of $A$.
\end{prop}



\begin{bem} If $A$ (and hence $B$) is symmetric, then $g=\sqrt {l}$ is
the minimizer in \eqref{ratefunction}, and we have
$I_A(\mu)=\|(-A)^{\frac 12}\sqrt{\mu}\|_2^2$. In this case the upper
bound simplifies to
\eq\label{rhoboundsymm}
    \rho_{ab}^{\ssup R}(l)
\leq {\rm e}^{-TI_A(\frac 1Tl)}
\Big(\prod_{x\in R\setminus \{a,b \}}\sqrt{\frac{T}{l_x}}\Big)\eta_{\sR}^{|R|-1}{\rm e}^{|R|[\chdb{1}+(4\eta_{\sR} T)^{-1}]}.
\en

\hfill$\Diamond$
\end{bem}

The proof of Proposition~\ref{lem-rhobound} makes use of three
lemmas that we will state and prove first.

\begin{lemma}\label{lem:monotone}
Let $\widetilde B\in [0,\infty)^{R\times R}$ be any \ch{matrix with nonnegative elements}, and let $Q\subseteq R$. Then
\begin{equation}\label{eq:monotone2}
    0
\leq
    \partial_{l}^{Q}\int_{[0,2\pi]^R}
    {\rm e}^{\sum_{x,y\in R} {\widetilde B}_{x,y}\sqrt{l_x}\sqrt{l_y}{\rm e}^{\i(\theta_{x}-\theta_{y})}}
    \prod_{x\in R}\frac{\d \theta_x}{2\pi}
\le
    \partial_{l}^{Q}
    {\rm e}^{\sum_{x,y\in R} {\widetilde B}_{x,y}\sqrt{l_x}\sqrt{l_y}},\qquad l\in(0,\infty)^R,
\end{equation}
where $\partial_{l}^{Q}=\prod_{x\in Q}\partial_{l_{x}}$.
\end{lemma}

\begin{proofsect}{Proof} Write ${\rm e}^{\sum_{x,y\in R}\dotsb }=
\prod_{x,y\in R}{\rm e}^{\dotsb}$ and expand the exponentials as power
series.  For $n= (n_{x,y})_{x,y\in R}\in\N_0^{R\times R}$, we write $n!=\prod_{x,y\in R} n_{x,y}!$. Then  we
obtain
\begin{equation}
\begin{aligned}
    \partial_{l}^{Q}& \int_{[0,2\pi]^R}
    {\rm e}^{\sum_{x,y\in R} {\widetilde B}_{x,y}\sqrt{l_x}\sqrt{l_y}{\rm e}^{\i(\theta_{x}-\theta_{y})}}
    \prod_{x\in R}\frac{\d \theta_x}{2\pi}
\\
&=
    \sum_{n\in \N_0^{R\times R}}\frac{1}{n!}\partial_{l}^{Q}\Big[\prod_{x,y\in R}
    \big(\widetilde B_{x,y}\sqrt{l_x}\sqrt{l_y}\big)^{n_{x,y}}
    \int_{[0,2\pi]^R}
    {\rm e}^{\i \sum_{x,y\in R} n_{x,y}(\theta_{x}-\theta_{y})}
    \prod_{x\in R}\frac{\d \theta_x}{2\pi}\Big].
\end{aligned}
\end{equation}
After rewriting the exponent in the integral on the right hand side using
$\sum_{x,y}n_{x,y}(\theta_{x}-\theta_{y})=\sum_{x}n_{x}\theta_{x}$,
where $n_{x}=\sum_{y} (n_{x,y}-n_{y,x})$, it is clear that the integral
equals one or zero. Hence, the lower bound in \eqref{eq:monotone2} is clear,
and the upper bound comes
from replacing the integral by one and a resummation over $n$.
\qed
\end{proofsect}

\begin{lemma}\label{lem:determinant-upper-bound} Fix any matrix $B\in\RR^{R\times R}$, let $a,b\in R$, and let $f\colon (0,\infty)^R\to\RR$ be any function with
nonnegative derivatives, i.e., $\partial_{l}^{Q}f
(l)\ge 0$ for all $Q\subseteq R$. Then
\begin{equation}\label{eq:determinant-upper-bound1}
    \big|
    \cofactor{a,b}^{\ssup{R}} (-B+\partial_{l})f
    \big|
\le
    \eta_{\sR}\prod_{x\in R\setminus \{a,b \}}\big(\eta_{\sR} +\partial_{l_{x}}\big)f,
\end{equation}
where $\eta_{\sR}$ is defined in \refeq{new-CRdef}.
\end{lemma}

\begin{proofsect}{Proof}
Recalling that the determinant is the (signed) volume subtended by the
rows, we can bound a determinant by the product of the lengths of the
rows. This is called the Hadamard bound and it applies to any real
square matrix.  Therefore, for $X\subseteq R$ and $a,b\in X$,
\[
    \big|\cofactor{ab}^{\ssup{X}}(-B)\big|
\le
    \prod_{x \in X\setminus\{\chdb{b}\}} \|B_{x}\|
\le
    \prod_{x \in X\setminus\{\chdb{b}\}} \eta_{\sR}=\eta_{\sR}^{|X|-1},
\]
where $B_{x}$ is \ch{the row $x$ of $B$} after eliminating
the \chdb{$a$}-th column, and $\|\cdot\|$ is the Euclidean length,
\chdb{which is bounded by $\eta_{\sR}$ because $\sum |a_{i}|^{2} \le
(\sum |a_{i}|)^{2}$.}
Also,
    \eq
    \cofactor{ab}^{\ssup{R}}(-B+\partial_l)f(l)
    =
    \sum_{\sigma\colon R\setminus\{\chdb{b}\}\to R\setminus\{\chdb{a}\}}
    \sign(\chdb{\hat{\sigma}})\prod_{x\in R\setminus\{a\}}\big(- B_{x,\sigma_x}
    +\delta_{x,\sigma_x}\partial_{l_x}\big)f(l),
    \en
where the sum over $\sigma$ is over all bijections
\chdb{$R\setminus\{b\}\to R\setminus\{a\}$}, \chdb{and where $\sign(\hat{\sigma})$
is the sign of the permutation $\hat{\sigma}\colon R\mapsto R$
obtained by letting $\hat{\sigma}_{x}=\sigma_{x}$ for $x\not =b$ and $\hat{\sigma}_{b}=a$.}
Expanding the product, we obtain
    \eq
    \begin{aligned}
    \cofactor{ab}^{\ssup{R}}(-B+\partial_l)f(l)
    &=
    \sum_{Q\subseteq R\setminus\{a,b\}} \ \sum_{\sigma\colon Q^{\rm c}
    \setminus\{\chdb{b}\}\to Q^{\rm c}\setminus\{\chdb{a}\}}\sign(\chdb{\hat{\sigma}})
    \Big(\prod_{x\in Q^{\rm c}\setminus\{\chdb{b}\}}\big(- B_{x,\sigma_x}\big)\Big)
    \Big(\prod_{x\in Q}\partial_{l_x}\Big)f(l)\\
    &=
    \sum_{Q\subseteq R\setminus\{a,b\}}\cofactor{ab}^{\ssup{Q^{\rm c}}}(-B)
    \Big(\prod_{x\in Q}\partial_{l_x}\Big)f(l),
    \end{aligned}
    \en
where we write $Q^{\rm c}=R\setminus Q$.  Take absolute values and bound the
cofactor using the Hadamard bound,
    \eq\lbeq{determinant-upper-bound2}
    \begin{aligned}
    \left|
    \cofactor{ab}^{\ssup{R}}(-B+\partial_l)f(l)
    \right|
    &\le\sum_{Q\subseteq R\setminus\{a,b\}}
    \eta_{\sR}^{|Q^{\rm c}\setminus \{\chdb{b}\}|}
    \Big(\prod_{x\in Q}\partial_{l_x}\Big)f(l)\\
    &=\eta_{\sR} \sum_{Q\subseteq R\setminus\{a,b\}}
    \Big(\prod_{x\in (R\setminus \{a,b \})\setminus Q}\eta_{\sR}\Big)
    \Big(\prod_{x\in Q}\partial_{l_x}\Big)f(l)\\
    &=
    \eta_{\sR} \prod_{x\in R\setminus \{a,b \}}\big(
    \eta_{\sR}+\partial_{l_{x}}
    \big)f(l).
    \end{aligned}
    \en
\qed
\end{proofsect}

\begin{lemma}\label{firstuppbound}
Fix any finite
subset $R$ of $\Lambda$, let $\widetilde B\in[0,\infty)^{R\times R}$ be any \ch{matrix with nonnegative elements}, and fix $a,b\in R$. Then, for any $T>0$ and any $l\in\Mcal^{+}_{\sT}$,
\begin{equation}\label{rhoboundfirst}
\begin{aligned}
    \cofactor{ab}^{\ssup{R}} &\big(-B+\partial_{l}\big)
    \int_{[0,2\pi]^R}
    {\rm e}^{\sum_{x,y\in R} \widetilde B_{x,y}\sqrt{l_x}\sqrt{l_y}{\rm e}^{\i(\theta_{x}-\theta_{y})}}
    \prod_{x\in R}\frac{\d \theta_x}{2\pi}\\
    &\le {\rm e}^{\sum_{x,y\in R} \widetilde B_{x,x}\sqrt{l_x}\sqrt{l_y}}
    \Big(\prod_{x\in R\setminus \{a,b \}}\sqrt{\frac{T}{l_x}}\Big)\eta_{\sR}^{|R|-1}{\rm e}^{[\chdb{\eta_{\sR}}^{-1}+(4\eta_{\sR}^2 T)^{-1}]\sum_{x,y\in R} \widetilde B_{x,y}},
    \end{aligned}
    \end{equation}
    where $\eta_{\sR}$ is defined in \refeq{new-CRdef}.
    \end{lemma}

\begin{proofsect}{Proof}
By Lemma~\ref{lem:determinant-upper-bound}
followed by Lemma~\ref{lem:monotone},
we obtain
    \eq \lbeq{new-upper-bound1}
    \mbox{l.h.s.~of \eqref{rhoboundfirst}}\le
    \eta_{\sR}\prod_{x\in R\setminus \{a,b \}}\big(\eta_{\sR} +\partial_{l_{x}}\big){\rm e}^{\sum_{x,y\in R} \widetilde B_{x,y}\sqrt{l_x}\sqrt{l_y}}.
    \en

Substitute  $t_{x}= \frac{\sqrt{l_{x}}}{\sqrt{T}}\in [0,1]$ and abbreviate $f (t)={\rm e}^{T\sum_{x,y\in R} \widetilde B_{x,y}t_{x}t_{y}}$. By the chain rule,
$\partial_{l_{x}}= \frac{1}{2T}\frac{1}{t_{x}}\partial_{t_{x}}$. Then
    \eq
    \begin{aligned}
    \mbox{l.h.s.~of \eqref{rhoboundfirst}}
    &\leq
    \eta_{\sR}^{|R|-1}
    \prod_{x\in R\setminus \{a,b \}}
    \Big(1 +\frac{1}{2\eta_{\sR} T}\frac{1}{t_{x}}\partial_{t_{x}}\Big)f(t)\\
    &\leq\eta_{\sR}^{|R|-1}\Big(\prod_{x\in R\setminus \{a,b \}}\frac 1{t_x}\Big)\prod_{x\in R\setminus \{a,b \}}
    \Big(1 +\frac{1}{2\eta_{\sR} T}\partial_{t_{x}}\Big)f(t),
    \end{aligned}
    \en
where we \ch{have} used that $t_x\leq 1$.
Since \emph{all} $t$ derivatives \ch{(not just the first order
derivatives)} of $f$ are nonnegative since
${\widetilde B}_{x,y}\geq 0$, we can add in some extra derivatives
and continue the bound with
    \eq
    \begin{aligned}
    \mbox{l.h.s.~of \eqref{rhoboundfirst}}&\le
    \eta_{\sR}^{|R|-1}\Big(\prod_{x\in R\setminus \{a,b \}}\frac 1{t_x}\Big)
    \prod_{x\in R\setminus \{a,b \}}
    \Big(\sum_{n=0}^\infty \frac{1}{n!} \frac{\partial^{n}_{t_{x}}}{(2\eta_{\sR} T)^{n}}\Big)f(t)\\
    &=
    \eta_{\sR}^{|R|-1}\Big(\prod_{x\in R\setminus \{a,b \}}\frac 1{t_x}\Big)
    \,f\Big(t+(2\eta_{\sR} T)^{-1}\1_\sR\Big),
    \end{aligned}
    \en
    where the last equation follows from Taylor's theorem, and $\1_\sR\colon R\to\{1\}$ is the constant function.

Recalling that $t_x\leq 1$, we may estimate
$$
\begin{aligned}
\frac 1T\log f\Big(t+(2\eta_{\sR} T)^{-1})\1_\sR\Big)
&=\sum_{x,y\in R}{\widetilde B}_{x,y}t_xt_y+\frac1{2\eta_{\sR} T}\sum_{x,y\in R}{\widetilde B}_{x,y}(t_x+t_y)+
 \frac1{(2\eta_{\sR} T)^2}\sum_{x,y\in R}{\widetilde B}_{x,y}
 \\
 &\leq \frac 1T\log f(t)+\frac 1T\Big[\frac1{\chdb{\eta}_{\sR} }+\frac1{4\eta_{\sR}^2T}\Big]\sum_{x,y\in R}{\widetilde B}_{x,y}.
\end{aligned}
$$
We conclude that
\begin{equation}
\mbox{l.h.s.~of \eqref{rhoboundfirst}}\le
\eta_{\sR}^{|R|-1}\Big(\prod_{x\in R\setminus \{a,b \}}\frac 1{t_x}\Big)\, f(t)\,{\rm e}^{[\chdb{\eta}_{\sR}^{-1}+(4\eta_{\sR}^2 T)^{-1}]\sum_{x,y\in R} {\widetilde B}_{x,y}}.
\end{equation}
Re-substituting $t_x=\ch{\sqrt{l_x/T}}$ and
$f(t)={\rm e}^{T\sum_{x,y\in R} {\widetilde B}_{x,y}t_{x}t_{y}}$, the lemma is proved.
\qed
\end{proofsect}

\begin{proofsect}{Proof of Proposition~\ref{lem-rhobound}}
Fix any $r\in(0,\infty)$ and recall the representation of the density $\rho^{\ssup R}_{ab}$ in \refeq{density3}. Now apply Lemma~\ref{firstuppbound} for $\widetilde B=(r_{x}B_{x,y}r_{y}^{-1})_{x,y\in R}$, to obtain
$$
  \rho^{\ssup R}_{ab}(l)\leq {\rm e}^{\sum_{x,y\in R} r_x\sqrt{l_x} A_{x,y}\sqrt{l_y}r_y^{-1}}
\Big(\prod_{x\in R\setminus \{a,b \}}\sqrt{\frac{T}{l_x}}\Big)\eta_{\sR}^{|R|-1}{\rm e}^{[\chdb{\eta}_{\sR}^{-1}+(4\eta_{\sR}^2 T)^{-1}]\sum_{x,y\in R} \widetilde B_{x,y}}.
    $$
Now we choose $r=\sqrt l/g$, where $g\in(0,\infty)^R$ is the minimiser in \eqref{ratefunction} for $\mu=\frac 1Tl$. This implies the bound in \eqref{rhobound}.
\qed\end{proofsect}

\subsection{Upper bounds in the LDP and in Varadhan's \chwk {l}emma}\label{sec-uppLDP}

In this section we specialize to Markov chains having a symmetric generator $A$ and give a simple upper bound for the left hand side of \refeq{upperLDP} and for
the expectation in \refeq{Varadhan}. Recall from the text below \eqref{ratefunction} that, in the present case of a symmetric generator, $I_A(\mu)=\|(-A)^{\frac 12}\sqrt \mu\|_2^2$ for any probability measure $\mu$ on $\Lambda$.

\begin{theorem}[Large deviation upper bounds for the local times]
\label{thm-DV}
Let the assumptions of Theorem~\ref{thm-locdens} be satisfied. Assume that $A$ is symmetric. Fix a finite subset $S$ of $\Lambda$. Then, for any $T\geq 1$ and any $a\in S$, with $\eta_{\sS}$ as in \refeq{new-CRdef}, \ch{the following bounds hold:}
\begin{itemize}
\item[{\rm (i)}]  \ch{For} every measurable $\Gamma\subseteq \Mcal_1(S)$,
    \eq\lbeq{uppprob}
    \log\mathbb{P}_a\big({\textstyle{\frac 1T}}\ell_{\sT}\in \Gamma, R_{\sT}\subseteq S\big)
    \leq -T\inf_{\mu\in \Gamma} \big\|(-A)^{\frac 12}\sqrt \mu\big\|_2^2+|S|\log \big(\eta_{\sS}\sqrt{8{\rm e}} T\big)+\log|S|+\frac {|S|}{4T}.
    \en

\item[{\rm (ii)}] For every measurable functional $F\colon\Mcal_1(S) \to \mathbb{R}$,
    \eq\label{expint}
    \log \mathbb{E}_a\big[{\rm e}^{T F(\frac 1T\ell_{\sT})}\1_{\{R_{\sT}\subseteq S\}}\big]
    \leq T\sup_{\mu\in \Mcal_1(S)}\Big[ F(\mu)-\big\|(-A)^{\frac 12}\sqrt \mu\big\|_2^2\Big]+|S|\log \big(\eta_{\sS}\sqrt{8{\rm e}} T\big)+\log|S|+\frac {|S|}{4T}.
    \en
\end{itemize}
\end{theorem}

Theorem~\ref{thm-DV} is a significant improvement over the standard estimates
known in large deviation theory. In fact, one standard technique to derive
upper bounds for the left hand side of \refeq{uppprob} is the use of the
exponential Chebyshev inequality and a compactness argument if $\Gamma$ is
assumed closed. One important ingredient there is a good control on the
logarithmic asymptotics of the expectation in \eqref{expint} for linear
functions $F$. This technique produces an error of order ${\rm e}^{o(T)}$,
which can in general not be controlled on a smaller scale.

The standard technique to derive improved bounds on the expectation in \eqref{expint}
for fixed $T$
is restricted to {\it linear\/} functions $F$, say $F(\cdot)=\langle V,\cdot\rangle$.
This technique goes via an eigenvalue expansion for the operator $A+V$ in the
set $S$ with zero boundary condition. The main steps are the use of the
Rayleigh-Ritz principle for the identification of the principal eigenvalue,
and Parseval's identity. This gives basically the same result as in
\eqref{expint}, but is strictly limited to linear functions $F$.

\vskip0.3cm

\begin{proofsect}{Proof of Theorem~\ref{thm-DV}} It is clear that (ii) follows from (i), hence we only prove (i).

According to Theorem~\ref{thm-locdens}, we may express the probability on the left hand side of \refeq{uppprob} as
    \eq\lbeq{start}
    \mathbb{P}_a\big({\textstyle{\frac 1T}}\ell_{\sT}\in \Gamma, R_{\sT}\subseteq S\big)
    =\sum_{b\in S}\sum_{R\subseteq S\colon a,b\in R}\int_{\Mcal^+_{\sT}(R)\cap
    \Gamma_{\sT,\sR}}\rho_{ab}^{\ssup{R}}(l)\,\sigma_\sT(\d l),
    \en
where $\Gamma_{\sT,\sR}=T\Gamma_{\sR}$, and $\Gamma_{\sR}$ is the set of the restrictions
of all the elements of $\Gamma$ to $R$.

We fix $\ch{a}, b\in S$ and $R\subseteq S$ with $a,b\in R$ and use the bound in Proposition~\ref{lem-rhobound}, more precisely, the one in \eqref{rhoboundsymm}.
Hence, for $l\in \Mcal^+_{\sT}(R)\cap \Gamma_{\sT,\sR}$, we obtain, after a substitution $l=T\mu$ in the exponent, that
    \eq
    \rho_{ab}^{\ssup R}(l)\leq
    {\rm e}^{-T\inf_{\mu\in \Gamma\colon \supp(\mu)\subseteq R}
    \|(-A)^{\frac 12}\sqrt \mu\|_2^2}\Big(\prod_{x\in R\setminus \{a,b \}}\sqrt{\frac{T}{l_x}}\Big)\eta_{\sR}^{|R|-1}{\rm e}^{|R|[\chdb{1} +(4\eta_{\sR} T)^{-1}]}.
    \en
Substituting this in \refeq{start} and integrating over $l\in \Mcal^+_{\sT}(R)$, we obtain
    \eq
    \begin{aligned}
    \mathbb{P}_a\big({\textstyle{\frac 1T}}
    \ell_{\sT}\in \Gamma, R_{\sT}\subseteq S\big)
    &\leq {\rm e}^{-T\inf_{\mu\in \Gamma}
    \|(-A)^{\frac 12}\sqrt \mu\|_2^2}
    \eta_{\sR}^{|R|-1}{\rm e}^{|R|[\chdb{1} +(4\eta_{\sR} T)^{-1}]}\\
    &\qquad\qquad\qquad\qquad\times\sum_{b\in S}
    \sum_{R\subseteq S\colon a,b\in R}
    \int_{\Mcal^+_{\sT}(R)}\prod_{x\in R\setminus\{a\}}\sqrt{\frac T{l_x}}\,\sigma_\sT(\d l)\\
    &\leq
    {\rm e}^{-T\inf_{\mu\in \Gamma}\|(-A)^{\frac 12}\sqrt \mu\|_2^2}\eta_{\sR}^{|R|-1}{\rm e}^{|S|[\chdb{1} +(4\eta_{\sR} T)^{-1}]}|S| {\sqrt 8}^{|S|}T^{|S|-1}.
    \end{aligned}
    \en
In the last integral, we \ch{have} eliminated $l_a=T-\sum_{y\in R\setminus \{a\}}l_y$,
\ch{have} extended the $(|R|-1)$ single integration areas to $(0,T)$ and used that
$\int_0^T l_x^{-\frac 12}\,\d l_x=\sqrt {2T}$. Now we use that $\eta_\sR$ is
increasing in $R$ and \ch{greater than or equal to}
one to arrive at \refeq{uppprob}. This completes the proof of (i).
\qed
\end{proofsect}

\subsection{Rescaled local times}\label{sec-rescaled}

As an application of Theorem~\ref{thm-DV}, we now consider continuous-time simple random walk restricted to a large $T$-dependent subset $\Lambda=\Lambda_{\sT}$ of $\Z^d$ increasing to $\Z^d$. We derive the sharp upper bound in the large deviation principle for its {\it rescaled\/} local times.  Assume, for some scale function  $T\mapsto{ \alpha_{\sT}}\in(0,\infty)$, that $\Lambda_{\sT}$ is equal to the box $[-R\alpha_{\sT},R\alpha_{\sT}]^d\cap\Z^d$, where the scale function $\alpha_{\sT}$ satisfies
    \eq\label{alphagrowth}
    1\ll\alpha_{\sT}\ll\Big(\frac T {\log{T}}\Big)^{\frac1{d+2}}\qquad\mbox{as }T\to\infty.
    \en
We introduce the rescaled version of the local times,
    $$
    L_{\sT}(x)=\frac{\alpha_{\sT}^d}T\ell_{\sT}\big(\lfloor \alpha_{\sT} x \rfloor\big), \qquad x\in \R^d.
    $$
Note that $L_{\sT}$ is a random step function on $\R^d$. In fact, it is a
random probability density on $\R^d$. Its support is contained in the cube
$[-R,R]^d$ if and only if the support of $\ell_\sT$ is contained in the box $[-R\alpha_{\sT},R\alpha_{\sT}]^d\cap\Z^d$.

It is known that, as $T\to\infty$, the family $(L_{\sT})_{T>0}$ satisfies a
large deviation principle under the sub-probability measures
$\P(\cdot\cap\{\supp(L_{\sT})\subseteq[-R,R]^d\})$ for any $R>0$. The speed
is $T\alpha_{\sT}^{-2}$, and the rate function is the energy \ch{functional}, i.e., the
map $g^2\mapsto \frac12\|\nabla g\|_2^2$, restricted to the set of squares
$g^2$ of $L^2$-normalized functions $g$ such that $g$ lies in $H^1(\R^d)$
and has its support in $[-R,R]^d$.
The topology is the one which is induced by all the test integrals of $g^2$
against
continuous and bounded functions. This large-deviation principle is proved
in \cite{GKS04} for
the discrete-time random walk, and the proof for continuous-time walks is
rather similar \ch{(see also \cite{HKM04}, where the proof of this fact
is sketched)}.
Hence, Varadhan's \chwk{l}emma yields precise logarithmic asymptotics for
all exponential functionals of $L_\sT$ that are bounded and continuous in
the above mentioned topology.

Note that this large deviations principle for $L_\sT$ is almost
the same as the one which is satisfied by the normalized Brownian
occupation times measures (see \cite{Ga77, DV75-83}), the main
difference being the speed (which is $T$ in \cite{Ga77, DV75-83}
instead of $T\alpha_{\sT}^{-2}$ here)
and the fact that $L_\sT$ does not take values in the set of
continuous functions $\R^d\to[0,\infty)$.

Here we want to point out that Theorem~\ref{thm-DV} yields a
new method to derive upper bounds for many exponential
functionals of $L_{\sT}$. For a cube $Q\ch{\subset \R^d}$, we denote by $M_1(Q)$
the set of all probability densities $Q\to[0,\infty)$.

\begin{theorem}\label{lem-rescloctim}Fix $R>0$, \ch{denote} $Q_{\sR}=[-R,R]^d$
and fix a measurable function
$F\colon M_1(Q_{\sR})\to\R$. Introduce
\eq
    \chi=\inf\Big\{\frac 12\|\nabla g\|_2^2-F(g^2)\colon g\in H^1(\R^d), \|g\|_2=1,\supp(g)\subseteq Q_{\sR}\big\}.
    \en
Then
    \eq\label{LDPscaledupper}
    \limsup_{T\to\infty}\frac{\alpha_{\sT}^2}
    T\log\E\Big[\exp\Big\{\frac T{\alpha_{\sT}^2}F(L_{\sT})\Big\}\1_{\{\supp(L_{\sT})\subseteq Q_{\sR}\}}\Big]
    \leq -\chi,
    \en
provided that
\begin{equation}\label{condition}
\liminf_{T\uparrow \infty}\inf_{\mu\in\Mcal_1(B_{R\alpha_{\sT}})}\Big( \alpha_{\sT}^2\frac 12
    \sum_{x\sim y}\Big(\sqrt{\mu(x)}-\sqrt{\mu(y)}\Big)^2-F\Big({\alpha_{\sT}^d}\mu\big(\lfloor \cdot\, \alpha_{\sT}\rfloor\big)\Big)\Big)
\geq \chi.
\end{equation}
\end{theorem}
\begin{proofsect}{Proof} Introduce
$$
F_{\sT}(\mu)=\frac 1{\alpha_{\sT}^2}F\Big({\alpha_{\sT}^d}\mu\big(\lfloor \cdot\, \alpha_{\sT}\rfloor\big)\Big),\qquad \mu\in\Mcal_1(\Z^d),
$$
then we have $\frac 1{\alpha_{\sT}^2}F(L_{\sT})=F_{\sT}(\frac 1T\ell_{\sT})$. Hence, Theorem~\ref{thm-DV}(ii) yields that
    $$
    \begin{aligned}
    \E\Big[\exp\Big\{\frac T{\alpha_{\sT}^2}F(L_{\sT})\Big\}\1_{\{\supp(L_{\sT})\subseteq Q_{\sR}\}}\Big]&=
    \E\Big[\exp\Big\{T \, F_{\sT}({\scriptstyle{\frac 1T}}\ell_{\sT})\Big\}\1_{\{\supp(\ell_{\sT})\subseteq Q_{R\alpha_{\sT}}\}}\Big]\\
    &\leq {\rm e}^{o(T\alpha_{\sT}^{-2})}{\rm e}^{-T  \chi_{\sT}},
    \end{aligned}
    $$
where
    $$
    \chi_{\sT}=\inf_{\mu\in\Mcal_1(Q_{R\alpha_{\sT}}\cap\Z^d)}\Big( \frac 12
    \sum_{x\sim y}\Big(\sqrt{\mu(x)}-\sqrt{\mu(y)}\Big)^2-F_{\sT}(\mu)\Big).
    $$
Here we used that the two error terms on the right hand side of \eqref{expint} are ${\rm e}^{o(T\alpha_{\sT}^{-2})}$ since $\eta_{\sS}\leq 2d$ for any $S\subseteq \Z^d$ and because of our growth assumption in \eqref{alphagrowth}. Now \eqref{LDPscaledupper} follows from \eqref{condition}.
\qed\end{proofsect}

\medskip

Theorem~\ref{lem-rescloctim} proved extremely useful in the study of the parabolic Anderson model in \cite{HKM04}. Indeed, it was crucial in that paper to find the precise upper bound of the left hand side of \eqref{LDPscaledupper} for the functional
$$
F(g^2)=\int_{Q_{\sR}} g^2(x)\log g^2(x)\,\d x,
$$
which has bad continuity properties in the topology in which the above mentioned large deviations principle holds. However, Theorem~\ref{lem-rescloctim} turned out to be applicable since the crucial prerequisite in \eqref{condition} had been earlier provided in \cite{GH99}. The main methods there were equicontinuity, uniform integrability and Arzela-Ascoli's theorem.

In the same paper \cite{HKM04}, also the functional
$$
F(g^2)=-\int_{Q_{\sR}}|g(x)|^{2\gamma}\,\d x,\qquad\mbox{with some }\gamma\in(0,1),
$$
was considered. This problem arose in the study of the parabolic Anderson model for another type of potential distribution which was earlier studied in \cite{BK01}. The prerequisite in \eqref{condition} was provided in \cite{HKM04} using techniques from Gamma-convergence; see \cite{AC04} for these techniques.

\section{Discussion}
\label{sec-dis}
In this section, we give some comments on the history of the problem
addressed in the present paper.

\subsection{Historical background}

The formulas in this paper have been motivated by the work of the
theoretical physicist J.M.~Luttinger \cite{Lut83} who gave a
(non-rigorous) asymptotic evaluation of certain path integrals.
Luttinger claimed that there is an asymptotic series
\[
    \E \big[{\rm e}^{-T F (\ell_{\sT}/T)}\big] \sim \sqrt{T}{\rm e}^{-c_{0}T}
    \Big(
    c_{1}+ \frac{c_{2}}{T} + \frac{c_{3}}{T^{2}} + \dotsb
    \Big)
\]
for Brownian local times.  He provided an algorithm to compute all the
coefficients.  He showed that his algorithm gives the Donsker-Varadhan
large deviations formula for $c_{0}$ and he explicitly computed the
central limit correction $c_{1}$.

In \cite{BrMu91} Brydges and Mu\~{n}oz{-}Maya used Luttinger's methods
to verify that his asymptotic expansion is valid to all orders for a
Markov process with symmetric generator and finite state space.  The
hypotheses are that $F$ is smooth and the variational principle that
gives the large deviations coefficient $c_{0}$ is non-degenerate.
\chdb{Luttinger implicitly relies on similar assumptions when he uses
the Feynman expansion for his functional integral}.

Thus there remains the open problem to prove that Luttinger's series
is asymptotic for more general state spaces, in particular, for
Brownian motion. As far as we know, the best progress to date is in
\cite{BoltDeu95} where compact state spaces were considered and the
asymptotics including the $c_{1}$ correction was verified.

Luttinger's paper used a calculus called \emph{Grassman integration}.
The background to this is that the Feynman-Kac formula provides a
probabilistic representation for the propagation of elementary
particles that satisfy ``Bose statistics''. To obtain a similar
representation for elementary particles that satisfy ``Fermi
statistics'' one is led in \cite{Bere87} to an analogue of integration
defined as a linear functional on a non-Abelian Grassman algebra in
place of the Abelian algebra of measurable functions: this is Grassman
integration. An important part of this line of thought concerns a case
where there is a relation called \emph{supersymmetry}.  This
background gives no hint that Grassman integrals are relevant for
ordinary Markov processes, but, nevertheless, Parisi and Sourlas
\cite{PaSo80} and McKane \cite{McK80} noted that random walk
expectations can be expressed in terms of the Grassman extension of
Gaussian integration.  Luttinger followed up on these papers by being
much more explicit and precise about the supersymmetric representation
in terms of Grassman integration and by deriving his series.

In \cite{LJ87} Le Jan pointed out that Grassman integration in this
context is actually just ordinary integration in the context of
differential forms. The differential forms are the non-Abelian algebra
and the standard definition of integration of differential forms
provides the linear functional.  Since integration over differential
forms is defined in terms of ordinary integration one can remove the
differential forms, as we have done in this paper, but this obscures
the underlying mechanism of supersymmetry.  The formalism with
differential forms is explained in \cite[page 551]{BrIm03b} where it
is used to study the Green's function of a self-repelling walk on a
hierarchical lattice. Two other applications of the same formalism are
the proof of the Matrix-Tree \chwk{t}heorem in \cite{Abde03} and a result on
self-avoiding trees given in \cite{BrIm03c}.

Luttinger found an instance of a relation between the local time of a
Markov process on a state space $E$ and the square of a {\it Gaussian
field\/} indexed by $E$.  The first appearance of such a relation was
given by Symanzik in \cite{Symanzik69}. His statement is that the
sum of the local times of an ensemble of Brownian loops is the square
of a Gaussian field. The references given above to Parisi-Sourlas,
McKane and Luttinger removed the need for an ensemble by bringing, in
its place, Grassman integration.  The paper of Symanzik was not
immediately rigorous because he claimed his result for Brownian motion
but it makes almost immediate sense for Markov processes on finite
state spaces only. Based on this work a rigorous relation between the
square of a Gaussian field and local time of a random walk on a
lattice was given by Brydges, Fr\"ohlich and Spencer in \cite{BFS82}.
Dynkin \cite{Dynk83,Dynk84b,Dynk84a} showed that the identities of
that paper can be extended to Brownian motion in one and two
dimensions.  In this form, the {\em Dynkin Isomorphism}, it became a useful
tool for studying local time of diffusions and much work has been done
by Rosen and Marcus in exploiting and extending these ideas, e.g., see
\cite{MarRos96, EKMRS00}.  The relation between the local time and
the square of a Gaussian field is concealed in this paper in
\refeq{baseq} which relates the local time $\ell$ to $l = |\phi|^{2}$
where $\phi$ is Gaussian.  This is more obvious when $\phi$ is
expressed as $\phi =u+\i v$ instead of in terms of polar coordinates
$\phi =\sqrt{l}{\rm e}^{\i\theta}$.

\subsection{\ch{Relation to the Ray-Knight \chwk{t}heorem}.}
\label{sec-RK}
\chwk{Our density formula in Theorem~\ref{thm-locdens} can also be used to prove a
version of the {\it Ray-Knight
theorem} for continous-time simple random walk on $\Z$. The
well-known Ray-Knight theorem
for one-dimensional Brownian motion, see \cite[Sections~XI.1-2]{RY91},
\cite[Sections~6.3-4]{KS91},
was originally proved in \cite{Kn63, Ra63}. It describes the Brownian
local times, observed at certain stopping times, as a homogeneous
Markov chain in the spatial parameter. Numerous deeper
investigations of this idea have been made, e.g., for general
symmetric Markov processes \cite{EKMRS00}, for diffusions
with fixed birth and death points on
planar cycle-free graphs \cite{EK93, EK96}, and on
the relations to Dynkin's isomorphism \cite{Shep85}, \cite{Eise94}.

The (time and space) discrete version of the Ray-Knight theorem, i.e., for simple random walk on $\Z$, was also introduced in \cite{Kn63}, however it turned out there that it is not the local times on the sites, but on the edges that enjoys a Markov property. This idea has been used or re-invented a couple of times, e.g., for applications to random walk in random environment \cite{KKS75}, to reinforced random walk \cite{Toth96}, and to random polymer measures \cite{GH93}.

In the present situation of continuous time and
discrete space, it turns out that the local times themselves form a
nice Markov chain. However, a proof appears to be missing.
In fact, up to our best knowledge, \cite{MS87} is the only
paper that provides (the outline of) a proof, but only for the special case
where the walk starts and ends in the same point. }

\Switchshort{We state the result here, but omit the proof. The proof will appear
in an extended version \cite{BHK06}. We first introduce some notation.}
For fixed $b\in \Z$, we denote
    \begin{equation}
    T_b^h = \inf\{t>0\colon \ell_t(b) > h\},\qquad h>0,
    \end{equation}
the right-continuous inverse of the map $t\mapsto\ell_t(b)$.
We denote by
    \eq
    I_0(h) = \sum_{i=0}^{\infty} \frac{h^{2i}}{2^i(i!)^2},
    \en
the modified Bessel function.

\begin{theorem}[Ray-Knight \chwk{t}heorem for continuous-time random walks]\label{thm-RK}
Let $\ell_\sT$ defined in \refeq{loctim} be the local times of
continuous-time simple random walk $(X_t)_{t>0}$ on $\Z$. Let $b\in\N$ and $h>0$.
\begin{enumerate}
\item[{\rm (i)}] Under $\P_0$, the process $\big(\ell_{T_b^h}(b-x)\big)_{x=0}^{b}$ is a time-homogeneous discrete-time Markov chain on $(0,\infty)$, starting at $h$,
with transition density given by
    \eq
    \lbeq{f}
    f(h_1,h_2) = {\rm e}^{-h_1-h_2} I_0\big(2\sqrt{h_1 h_2}\big),\qquad h_1,h_2\in (0,\infty).
    \en

\item[{\rm (ii)}] Under $\P_0$, the processes $\big(\ell_{T_b^h}(b+x)\big)_{x\in\N_0}$ and $\big(\ell_{T_b^h}(-x)\big)_{x\in\N_0}$
are time-homogeneous discrete-time Markov chains on $[0,\infty)$ with transition probabilities given by
    \eq
    \lbeq{fstar}
    P^\star (h_1,\d h_2) = {\rm e}^{-h_1}\delta_{0}(\d h_2)+{\rm e}^{-h_1-h_2} \sqrt{\frac{h_1}{h_2}} I_0'\big(2\sqrt{h_1 h_2}\big)\,\d h_2,\qquad h_1,h_2\in [0,\infty).
    \en
    \item[{\rm (iii)}] The three Markov chains in (i) and (ii) are independent.
\end{enumerate}
\end{theorem}
\vskip0.3cm

\noindent
We note that Theorem \ref{thm-RK}(ii) and an outline of its proof can
be found in \cite[(3.1-2)]{MS87}. This proof uses an embedding of
the random walk into a Brownian motion and the Brownian Ray-Knight theorem; we
expect that Theorem \ref{thm-RK}(i) and (iii) can also be proved along these lines.
\Switchshort{In the extended version \cite{BHK06},}
\Switchlong{In Appendix \ref{sec-app},} using the density formula of
Theorem~\ref{thm-locdens}, we provide a proof of Theorem \ref{thm-RK} that is independent
of the Brownian Ray-Knight theorem. This opens up the possibility of producing a
new proof of this theorem, via a diffusion approximation of the
Markov chains having the transition densities
in \refeq{f} and \refeq{fstar}. Furthermore, we emphasize that our proof can also
be adapted to continous-time random walks on cycle-free graphs and has some potential
to be extended to more general graphs. Theorem~\ref{thm-locdens} contains
far-ranging generalisations of the Ray-Knight idea, which are to be studied in future.

\subsection*{Acknowledgment}

{ DB would like to thank the Natural
Sciences and Engineering Research Council of Canada for supporting his
research.} The work of RvdH was supported in part by
Netherlands Organisation for Scientific Research (NWO). WK would like
to thank the German Science Foundation for awarding a Heisenberg grant
(realized in 2003/04). This project was initiated during an extensive visit of
RvdH to the University of British Columbia, Vancouver, Canada.

\Switchlong{\renewcommand{\thesection}{\Alph{section}}
\setcounter{section}{0}

\numberwithin{equation}{section}
\numberwithin{theorem}{section}
\section{Appendix: Proof of Theorem~\ref{thm-RK}}
\label{sec-app}
In this section, we prove Theorem~\ref{thm-RK}. (Recall the discussion in Section \ref{sec-RK}.)
To prove Theorem~\ref{thm-RK}, we will need the following proposition, which is
of independent interest. Recall that the matrix $B$ is the off-diagonal part of the generator  $A$ of the Markov chain. Let
\eq\lbeq{define-f}
    g_{x,y}(t)
=
    \int_{[0,2\pi]}
    {\rm e}^{t (
    B_{x,y} \, {\rm e}^{\i\theta} +
    B_{y,x} \, {\rm e}^{-\i\theta}
    )}\,\, \frac{\d\theta}{2\pi},\qquad t>0,x,y\in\Z.
\en

The following does not assume that we are dealing with simple random walk,
but holds for any continuous-time nearest-neighbor Markov chain on $\Z$. That is, the generator $A$ is a {\it tridiagonal\/} matrix, which means that $A_{x,y}=0$ for any $x,y$ satisfying $|x-y|\geq 2$.

\begin{prop}\label{prop:DimOneRho} Assume that the conservative  generator $A=(A_{x,y})_{x,y\in\Z}$ of the Markov chain is a tridiagonal matrix. Let $R\subset \Z$ be a finite interval and let $a, b \in R$ with $a\leq b$. Then
    \begin{equation}\label{eq:DimOneRho1}
    \rho^{\ssup R}_{ab}(l)
    =
    {\rm e}^{\sum_{x\in R} A_{x,x}l_x}
    \Big[
    \prod_{\heap{x<a}{x\in R}} \partial_{l_{x}} g_{x,y}(\sqrt{l_{x}l_{y}})
    \Big]
    \Big[
    \prod_{a\le x < b} g_{x,y}(\sqrt{l_{x}l_{y}})
    \Big]\Big[
    \prod_{\heap{y>b}{y\in R}} \partial_{l_{y}} g_{x,y}(\sqrt{l_{x}l_{y}})
    \Big],
    \end{equation}
where $y=x+1$.
\end{prop}

Proposition \ref{prop:DimOneRho} says that, for any nearest-neighbor Markov chain, the distribution of the sequence of local times possesses a product structure and may be divided into the piece between starting and ending point and the two boundary pieces. All the three pieces have Markovian structure, with explicit identification of the transition probability function, which is in general inhomogeneous.

The proof of Proposition \ref{prop:DimOneRho} makes use of the following two lemmas. The first one is purely algebraic and shows that the determinant of a tridiagonal matrix, after erasing one row and one column, naturally decomposes in a product of three parts.

\begin{lemma}\label{lem:tridiag} Let $R\subset\Z$ be a finite interval and let $a, b \in R$ with $a\leq b$.  Let $M$ be an
$R\times R$ tridiagonal matrix. Then
    \eq \lbeq{tridiag1}
    \cofactor{ab}(M)
    =
    \cofactor{aa}^{\ssup{\le a}}(M) \
    \Big(
    \prod_{a \le i <b} M_{i+1,i}\Big) \
    \cofactor{bb}^{\ssup{\ge b}}(M),
    \en
where the superscripts $\le a$ and $\ge b$ denote the sets $R^{\ssup{\le a}}=\{i\in
R\colon i\le a \}$ and $R^{\ssup{\geq b}}=\{i\in R\colon i\ge b\}$ respectively. If the product is
over the empty set or if the cofactors apply to empty matrices, then
the corresponding factor is set equal to $1$.
\end{lemma}

\begin{proofsect}{Proof} Let $M^{\ssup{<a}}
= (M_{i,j})_{i,j<a}$ and likewise for other inequalities as
superscripts.  Let $\widehat{M}$ denote the matrix obtained by removing
row $a$ and column $b$ from $M$ so that $\widehat{M}_{i,j} = M_{f (i),g
(j)}$, where $f (i)= i$ for $i<a$ and $f (i)=i+1$ for $i\ge
a$. Likewise $g (j) = j$ for $j<b$ and $g (j)=j+1$ for $j\ge b$.
{ In the proof, we distinguish two separate cases, depending
on whether $a=b$ or $a<b$.}

\emph{Case } $a=b$.  Then it is easy to see that
$\widehat{M}=M^{\ssup{<a}}\oplus M^{\ssup{>a}}$ is block-diagonal.
Since $\detrm(
M^{\ssup{< a}})=\cofactor{aa}^{\ssup{\le a}}(M)$ and $\detrm
(M^{\ssup{>a}})=\cofactor{bb}^{\ssup{\ge a}}(M)$, whereas the middle term on the
right of \refeq{tridiag1} is an empty product, the result
\refeq{tridiag1} follows immediately.

\emph{Case } $a<b$. The $R^{\ssup{<a}}\times R^{\ssup{<a}}$-submatrix of $\widehat{M}$ (the upper left corner)
is equal to $M^{\ssup{<a}}$, and right of this block there
are throughout zeros in $\widehat{M}$,
except for the last row (the row indexed by $a$).
Analogously, the $R^{\ssup{\geq b}}\times R^{\ssup{\geq b}}$-submatrix
of $\widehat{M}$ (the lower right corner)
is equal to $\widehat M^{\ssup{\geq b}}=M^{\ssup{> b}}$,
and there are throughout zeros above that block in $\widehat{M}$,
except for the first column (the column indexed by $b$). The intermediate block
$(\widehat{M}_{i,j})_{a\leq i,j<b}$
is in upper triangular form since, for $i=a,\dotsc ,b-1$ and $j>i$,
we have $\widehat{M}_{i,j}= M_{f
(i),g(j)}=M_{i+1,j}=0$. Furthermore, the only non-zero entries below the diagonal of
$\widehat{M}$ are in the diagonal that is next to the main diagonal.

We are going to calculate the determinant of $\widehat{M}$ by applying
linear row and column  operations that transform $\widehat{M}$ into upper triangular form.
For $i=\min R,\dots,a-2$, we add a suitable multiple of the $i$-th row to
the $(i+1)$-st row in order that the upper left corner is turned into a
upper triangular matrix. Note that these operations do not affect any
entry outside this corner.
Furthermore, for $j=\max R,\max R-1,\dots,b+1,$ we add a suitable
multiple of the $j$-th column to the $(j-1)$-st column in order that the
lower right corner is turned into an
upper triangular matrix. Note that these operations do not affect any
entry outside this corner.

The resulting $R\times R$-matrix is in upper triangular form, and,
since row additions do not change the determinant, its determinant is
equal to
$\det(\widehat M)=\det_{ab}(M)$. This determinant is equal to the
product of the three determinants of the left upper corner (which is $\detrm
(M^{\ssup{<a}})=\cofactor{aa}^{\ssup{\le a}}(M)$), the right lower corner
(which is $\detrm
(M^{\ssup{> b}})=\cofactor{bb}^{\ssup{\ge b}}(M)$) and the product of the
diagonal entries of the piece inbetween, which is $\prod_{i=a}^{b-1} \widehat M_{i,i}
=\prod_{a \le i <b} M_{i+1,i}$. This ends the proof.
\qed

\end{proofsect}

Now we state and prove the second lemma that will be used in the proof of Proposition~\ref{prop:DimOneRho}.

\begin{lemma}\label{lem:cofactorzero} Let $V$ be a finite interval in $\Z$ with $ |V|\geq 2$ and let $c\in V$ be the smallest
or the largest state. 
Then, for every $l\in \R^{V}_{+}$,
\begin{eqnarray}
    \cofactor{cc}^{\ssup{V}} \big(-B+\partial_{l}\big)
    \Big[
    \prod_{x,y\in V\setminus \{c\}}
    g_{x,y}(\sqrt{l_{x}l_{y}})
    \Big]
&=&
    0,\label{eq:cofactorzero1}\\
    \cofactor{cc}^{\ssup{V}} \big(-B+\partial_{l}\big)
    \Big[
    \prod_{x,y\in V} g_{x,y}(\sqrt{l_{x}l_{y}})
    \Big]
&=&
    \prod_{x,y\in V} \partial_{l_{z}}g_{x,y}(\sqrt{l_{x}l_{y}}),
    \label{eq:cofactorzero2}
\end{eqnarray}
where $y=x+1$ and where $z=x$ if $c$ is the largest state and $z=y$ if
$c$ is the smallest state in $V$.

\end{lemma}
\vskip0.3cm

\begin{proofsect}{Proof}
We first prove (\ref{eq:cofactorzero1}). In Theorem~\ref{thm-locdens}, choose $F (l)=f
(l_{c})$ to depend only on $l_{c}$ and integrate over $T$, to obtain that
\[
    \int_{0}^{\infty} \d T \,\,
    \E_c\big[f (\ell_{\sT}(c))\1_{\{X_{\sT}=c\}}\1_{\{R_{\sT}=V\}}\big]
=
    \int \d l_{c}\,\,f(l_{c})
    \int \rho^{\ssup V}_{cc}(l) \,\d^{V\setminus \{c \}}l.
\]
In the left hand side the indicator functions of the events
$\{R_{\sT}=V\}$ and $\{ {X_{\sss 0}}=X_{\sT}=c\}$ (recall that $|V|\geq 2$)
require the Markov chain to visit $c$
at least two times.  Conditioning on the number of visits to $c$ to be
$n_{c}$, the amount of time $\ell_{\sT}(c)$ spent at $c$ is the sum of
$n_{c}\ge 2$ exponential random variables with parameter
$A_{c,c}$. The density at $0$ of a sum of two or more exponential
random variables is zero. Therefore, for $l_{c}=0$,
\[
    \int \rho^{\ssup V}_{cc}(l) \,\d^{V\setminus \{c \}}l
=
    0.
\]
It can be seen from \refeq{density2} that $l_c\mapsto
\rho^{\ssup V}_{cc}(l)$ is continuous at $l_{c}=0$. Therefore,
$\rho^{\ssup V}_{cc}(l) =0$ at $l_{c}=0$. Also by \refeq{density2},
for $l_{c}=0$,
\begin{equation*}
    0
=
    \cofactor{cc}^{\ssup{V}} \big(-B+\partial_{l}\big)
    \Big[\prod_{x,y \in V} g_{x,y}(\sqrt{l_{x}l_{y}})\Big].
\end{equation*}
In $\cofactor{cc}^{\ssup{V}} \big(-B+\partial_{l}\big)$ no derivative
with respect to  $l_c$ appears. Also $g_{x,y}(0)=1$. Therefore, for $l_c=0$, we
may replace $V$ by $V\setminus \{c \}$ under the product sign. But this is (\ref{eq:cofactorzero1}).

Now we prove (\ref{eq:cofactorzero2}). The proof is by induction on the number of elements,
$n$, in $V$. For $n=2$ the statement is easily checked. This
initialises the induction, and we are left to advance it. We give the
argument only for the case where $c$ is the smallest vertex in $V$;
the other case is essentially the same. We split the product into
$$
\prod_{x,y \in V} g_{x,y}(\sqrt{l_{x}l_{y}})=
     g_{c,c+1}(\sqrt{l_{c}l_{c+1}})
     \prod_{x,y \in V\setminus \{c\}}g_{x,y}(\sqrt{l_{x}l_{y}})
$$
apply the differential operator $\cofactor{cc}^{\ssup{V}} (-B+\partial_{l})$
and use the product rule of differentiation for the derivative with
respect to $\partial_{l_{c+1}}$.  Since $\cofactor{cc}^{\ssup{V}}=\detrm^{\ssup{V\setminus
\{c \}}}$ is linear in the first row we therefore have
    \eqalign
    \lbeq{aimRKF1c}
    { \text{ l.h.s.~ of } \eqref{eq:cofactorzero2}}
    &=\big(\partial_{l_{c+1}}g_{c,c+1}(\sqrt{l_{c}l_{c+1}})\big)
     \Big(\cofactor{c+1\,c+1}^{\ssup{V\setminus\{c\}}} \big(-B+\partial_{l}\big)
     \prod_{x,y \in V\setminus \{c\}}g_{x,y}(\sqrt{l_{x}l_{y}})
     \Big)\nonumber\\
    &\qquad \quad+g_{c,c+1}(\sqrt{l_{c}l_{c+1}})
     \Big(
     \cofactor{cc}^{\ssup{V}} \big(-B+\partial_{l}\big)
     \prod_{x,y \in V\setminus \{c\}} g_{x,y}(\sqrt{l_{x}l_{y}})
     \Big). \nonumber
    \enalign
By the induction hypothesis, the first term is equal to
    \eq
    \big(\partial_{l_{c+1}}g_{c,c+1}(\sqrt{l_{c}l_{c+1}})\big)\prod_{x,y \in V\setminus \{c\}}
    \partial_{l_{y}}g_{x,y}(\sqrt{l_{x}l_{y}})
    =\prod_{x,y \in V} \partial_{l_{y}}g_{x,y}(\sqrt{l_{x}l_{y}}).
    \en
On the other hand, by (\ref{eq:cofactorzero1}), the second term is equal to zero.
This advances the induction, and, thus, completes the proof of (\ref{eq:cofactorzero2}).
\qed
\end{proofsect}

\begin{proofsect}{Proof of Proposition~\ref{prop:DimOneRho}}
Compare \refeq{define-f} and \refeq{density2} to see that
\begin{equation}
\lbeq{integralstheta}
\rho^{\ssup R}_{ab}(l)={\rm e}^{\sum_{x\in R} A_{x,x}l_x}
    \cofactor{ab}^{\ssup {R}}
    \big(-B+\partial_{l}\big)
    \prod_{x,y \in R} g_{x,y}(\sqrt{l_x l_{y}}),
\end{equation}
where we used our convention $y=x+1$. By Lemma~\ref{lem:tridiag} we
rewrite this as
    \begin{align}
    \lbeq{FRformula}
    \rho^{\ssup R}_{ab}(l)
    &=
    {\rm e}^{\sum_{x\in R} A_{x,x}l_x}
    \Big[
    \cofactor{aa}^{\ssup{\le a}} \big(-B+\partial_{l}\big)
    \prod_{x<a} g_{x,y}(\sqrt{l_{x}l_{y}})
    \Big]\\
    &\qquad\qquad \times
    \Big[
    \prod_{a\le x<b} g_{x,y}(\sqrt{l_{x}l_{y}})
    \Big]\Big[
    \cofactor{bb}^{\ssup{\ge b}} \big(-B+\partial_{l}\big)
    \prod_{y>b} g_{x,y}(\sqrt{l_{x}l_{y}})
    \Big].
    \nonumber
    \end{align}
By Lemma~\ref{lem:cofactorzero}, this equals the right-hand side of
(\ref{eq:DimOneRho1}).
\qed
\end{proofsect}

\begin{proofsect}{Proof of Theorem~\ref{thm-RK}}
We specialize Proposition~\ref{prop:DimOneRho} to simple random walk, whose generator $A=(A_{x,y})_{x,y\in\Z}$ is the Laplace operator, i.e., the tridiagonal matrix satisfying $A_{x,y}=1$ for $|x-y|=1$ and $A_{x,x}=-2$ for $x\in\Z$.
By \refeq{define-f} and since $B_{x,y}=B_{y,x}=1$,
    \eq
    g_{x,y}(t)
    =\int_{[0,2\pi]}{\rm e}^{2t \cos \theta }
     \,\, \frac{\d\theta}{2\pi}
    =I_0(2t).
    \en

Fix a finite interval $R\subseteq\Z$ containing $a,b$. Proposition~\ref{prop:DimOneRho} gives that, using that $A_{x,x}=-2$ and using the convention $l_{1+\max R}=l_{-1+\min R}=0$,
    \begin{equation}\lbeq{rhodl}
    \begin{aligned}
    \rho^{\ssup R}_{ab}(l)\,\d^Rl
    &=
    {\rm e}^{-l_{\min R}}\Big[\prod_{\heap{x<a}{x\in R}} {\rm e}^{-(l_x+l_y)} \partial_{l_{x}} I_0(2\sqrt{l_{x}l_{y}})\,\d l_x
    \Big]
    \Big[
    \prod_{a\le x < b} {\rm e}^{-(l_x+l_y)} I_0(2\sqrt{l_{x}l_{y}})\,\d l_x
    \Big]\\
    &\qquad\qquad\times\Big[\prod_{\heap{y>b}{y\in R}} {\rm e}^{-(l_x+l_y)}\partial_{l_{y}} I_0(2\sqrt{l_{x}l_{y}})\,\d l_y
    \Big]{\rm e}^{-l_{\max R}}\\
    &=\Big[\prod_{\heap{x\leq a}{x\in R}} P^{\star}(l_x,\d l_{x-1})\Big]
    \Big[
    \prod_{a\le x < b} f(l_{x}, l_{x+1})\,\d l_{x}
    \Big]\Big[\prod_{\heap{x\geq b}{x\in R}}P^{\star}(l_x,\d l_{x+1})
    \Big].
    \end{aligned}
    \end{equation}
This is the fixed-time equivalent of Theorem~\ref{thm-RK}.
To go to the stopping time $T_b^h$, we claim that, for any measurable set $C\subseteq(0,\infty)^{R\setminus \{b\}}$,
    \eq
    \lbeq{stoppingtime}
    \begin{aligned}
    {\mathbb P}\big(&(\ell_{T_b^h}(x))_{x\in R\setminus \{b\}}\in C,\supp(\ell_{T_b^h})=R, T_b^h\in \d T)\,\d h\\
    &={\mathbb P}\big((\ell_\sT(x))_{x\in R\setminus \{b\}}\in C,\supp(\ell_\sT)=R, X_{\sT}=b, \ell_{\sT}(b)\in \d h\big)\,\d T.
    \end{aligned}
    \en
Indeed, \refeq{stoppingtime} implies, for any bounded measurable functional $F\colon (0,\infty)^R\to\R$ and any bounded measurable
function $f\colon(0,\infty)\to\R$, that
    \eq
    \int_{0}^{\infty} f(h){\mathbb E}_a\big[F(\ell_{T_b^h})\1_{\{R_{T_b^h}=R\}}\big]\d h
    =\int_{0}^{\infty} {\mathbb E}_a\big[F(\ell_{\sT})\1_{\{R_\sT=R\}}\1_{\{X_{\sT}=b\}}f(\ell_{\sT}(b))\big]\,
    \d T.
    \en
By Theorem~\ref{thm-locdens} and \eqref{eq:dsigma}, the right hand side is equal to $\int_{(0,\infty)^R}F(l)f(l_b)\rho_{ab}^{\ssup R}(l)\,\d^R l$. Equation \refeq{rhodl}
gives that, on the event $\{\supp(\ell_{T_b^h})=R\}$, the distribution of $\ell_{T_b^h}$
is the one that is claimed in Theorem \ref{thm-RK}. Hence, Theorem~\ref{thm-RK} follows from \refeq{stoppingtime}.

We now prove \refeq{stoppingtime}. It is sufficient to prove that, for every $h>0$,
\eq
    \lbeq{stoppingtimeeps}
    \begin{aligned}
    {\mathbb P}\Big(&(\ell_{T_b^h}(x))_{x\in R\setminus \{b\}}\in C, \ch{\supp(\ell_{T_b^h})=R,} T_b^h\in (T,T+\eps)\Big)\\
    &={\mathbb P}\Big((\ell_{\sT}(x))_{x\in R\setminus \{b\}}\in C, \ch{\supp(\ell_{\sT})=R,} X_{\sT}=b, \ell_{\sT}(b)\in (h-\eps,h)\Big)+o(\eps),\qquad\mbox{as }\eps\downarrow 0.
    \end{aligned}
    \en
Equation \refeq{stoppingtimeeps} is a consequence of
    \eqalign
    \lbeq{stoppingtimeepsb}
    &{\mathbb P}(\ch{(\ell_{T_b^h}(x))_{x\in R\setminus \{b\}}\in C, \supp(\ell_{T_b^h})=R,}T_b^h\in (T,T+\eps))\nonumber\\
    &\quad ={\mathbb P}\big(\ch{(\ell_{T_b^h}(x))_{x\in R\setminus \{b\}}\in C, \supp(\ell_{T_b^h})=R,}
    X_{t}=b\, \forall t\in (T,T+\eps),\ell_{\sT}(b)\in (h-\eps,h)\big)+o(\eps).
    \enalign
Indeed, \refeq{stoppingtimeeps} follows from \refeq{stoppingtimeepsb} since,
\rm{conditionally} on the event $\{T_b^h\in (T,T+\eps)\}$, with high probability,
the random walker spends all the time in the interval $(T,T+\eps)$ in $b$ and
therefore does not change the local times in any other point than $b$ during
that time.

To prove \refeq{stoppingtimeepsb}, note that, on $\{T_b^h\in (T,T+\eps)\}$, we have $\ell_{\sT}(b) \leq  h < \ell_{\sss T+\eps}(b)$. Since also $\ell_{t+s}(b)\leq \ell_t(b)+s$ for any $t,s\geq 0$,
we also have that $\ell_{\sT}(b)\in (h-\eps,h]$. This shows that $\{T_b^h\in (T,T+\eps)\}\subseteq\{\ell_{\sT}(b)\in (h-\eps,h]\}$.
Furthermore, on $\{T_b^h\in (T,T+\eps)\}$, there is a $t\in (T,T+\eps)$ such that $X_t=b$.
On this event, the probability of the event $\{\exists s\in (T,T+\eps)\colon X_s\neq b\}$ is $\leq\Ocal(\eps^2)$, since $\ell_{\sT}(b)\in (h-\eps,h]$ on this event and at least one step happens during the time interval $(T,T+\eps)$.
Since $\ell_{\sT}(b)$ has a density, this first event has probability $\leq \Ocal(\eps)$ the second has probability $\leq \Ocal(\eps)$, and by the Markov property at time $T$, the intersection of these events has probability $\leq \Ocal(\eps^2)$. This ends the proof of \refeq{stoppingtimeepsb}.
\qed
\end{proofsect}
}


\begin{thebibliography}{WWW98}

\bibitem[Abd03]{Abde03}
\chwk{{\sc A.~Abdesselam,}
\newblock {Grassmann-Berezin calculus and theorems of the matrix-tree type}.
\newblock {\it Adv.~Appl.~Math.} {\bf 33:1}, 51--70 (2004).}

\smallskip

\bibitem[AC04]{AC04}
{\sc R.~Alicandro} and {\sc M.~Cicalese,}
\newblock {A general integral representation result for continuum
limits of discrete energies with superlinear growth.}
\newblock {\it SIAM J.\ Math.\ Anal.} {\bf 36}, 1--37 (2004).

\smallskip

\bibitem[Ber87]{Bere87}
{\sc F.A.~Berezin,}
\newblock {\it Introduction to Superanalysis}.
\newblock Mathematical Physics and Applied Mathematics {\bf 9}, D.~Reidel Publishing, Dordrecht
(1987).

\smallskip

\bibitem[BK01]{BK01}
{\sc M.~Biskup} and {\sc W.~K\"onig,}
\newblock {Long-time tails in the parabolic Anderson model
with bounded potential}.
\newblock {\em Ann. Probab.} {\bf 29:2}, 636-682 (2001).

\smallskip

\bibitem[BDT95]{BoltDeu95}
{\sc E.~Bolthausen, J.-D.~Deuschel}, and {\sc Y.~Tamura},
\newblock Laplace approximations for large deviations of nonreversible {M}arkov
  processes. {T}he nondegenerate case.
\newblock {\em Ann. Probab.} {\bf 23:1}, 236--267 (1995).

\smallskip

\bibitem[BFS82]{BFS82}
{\sc D.C. Brydges, J.~Fr\"{o}hlich,} and {\sc T.~Spencer,}
\newblock The random walk representation of classical spin systems and
  correlation inequalities.
\newblock {\em Commun. Math. Phys.} {\bf 83}, 123--150 (1982).

\smallskip

\Switchshort{
\bibitem[BHK06]{BHK06}
{\sc D.C. Brydges, R.~van der Hofstad,} and {\sc W.~K\"onig,}
\newblock Joint density for the local times of continuous-time Markov chains:
Extended version.
\newblock {\tt http:arXiv...}
}

\smallskip

\bibitem[BI03a]{BrIm03c}
{\sc D.C. Brydges} and {\sc J.Z. Imbrie,}
\newblock Branched polymers and dimensional reduction.
\newblock {\em Annals of Mathematics} {\bf 158}, 1019­--1039 (2003).

\smallskip

\bibitem[BI03b]{BrIm03b}
{\sc D.C. Brydges} and {\sc J.Z. Imbrie,}
\newblock Green's function for a hierarchical self-avoiding walk in four
  dimensions.
\newblock {\em Commun. Math. Phys.} {\bf 239:3}, 549--584 (2003).

\smallskip

\bibitem[BM91]{BrMu91}
{\sc D.C.~Brydges} and {\sc I.~Mu\~{n}oz{-}Maya},
\newblock An application of {B}erezin integration to large deviations.
\newblock {\em Jour. Theor. Probab.} {\bf 4}, 371--389 (1991).

\smallskip

\bibitem[DZ98]{DZ98}
{\sc A.~Dembo} and {\sc O.~Zeitouni,}
\newblock \textit{Large Deviations Techniques and Applications},
\newblock 2$^{\rm nd}$ edition. Springer, New York (1998).

\smallskip

\bibitem[DV75-83]{DV75-83}
{\sc M.D.~Donsker} und {\sc S.R.S.~Varadhan},
\newblock Asymptotic evaluation of certain Markov process expectations for large time, I--IV,
\newblock {\it Comm.~Pure Appl.~Math.} {\bf 28}, 1--47, 279--301 (1975), {\bf 29}, 389--461 (1979), {\bf 36}, 183--212 (1983).



\smallskip



\bibitem[Dyn83]{Dynk83}
{\sc E.M.~Dynkin},
\newblock {Gaussian and non-Gaussian random fields associated with Markov processes.}
\newblock {\it J.\ Funct.\ Anal.} {\bf 55}, 344--376 (1983).


\smallskip

\bibitem[Dyn84a]{Dynk84a}
{\sc E.M.~Dynkin},
\newblock {Local times and quantum fields.}
\newblock {\it Seminar on Stochastic Processes}, 1983 (Gainesville, Fla., 1983), 69--83,
{\it Progr.\ Probab.\ Statist.}, {\bf 7} Birkh\"auser Boston, Boston, MA (1984).


\smallskip

\bibitem[Dyn84b]{Dynk84b}
{\sc E.M.~Dynkin},
\newblock {Polynomials of the occupation field and related random fields.}
\newblock {\it J.\ Funct.\ Anal.} {\bf 58}, 20--52 (1984).

\smallskip

\bibitem[Eise94]{Eise94}
{\sc N.~Eisenbaum},
\newblock Dynkin's isomorphism theorem and the Ray-Knight
theorems.
\newblock {\em Probab.~Theory Relat.~Fields} {\bf 99}, 321--335 (1994).


\smallskip

\bibitem[EK93]{EK93}
{\sc N.~Eisenbaum} and {\sc H.~Kaspi},
\newblock A necessary and sufficient condition for the Markov
properties of the local time process.
\newblock {\em Ann. Probab.} {\bf 21:3}, 1591--1598 (1993).


\smallskip

\bibitem[EK96]{EK96}
{\sc N.~Eisenbaum} and {\sc H.~Kaspi},
\newblock On the Markov
property of local time for Markov processes on general graphs.
\newblock {\em Stoch.~Proc.~Appl.} {\bf 64}, 153--172 (1996).


\smallskip

\bibitem[EKMRS00]{EKMRS00}
{\sc N.~Eisenbaum, H.~Kaspi, M.~B. Marcus, J.~Rosen},
and {\sc Z.~Shi},
\newblock A {R}ay-{K}night theorem for symmetric {M}arkov processes.
\newblock {\em Ann.~Probab.} {\bf 28:4}, 1781--1796 (2000).


\smallskip

\bibitem[GKS05]{GKS04}
{\sc N.~Gantert, W.~K\"onig} and {\sc Z.~Shi},
\newblock Annealed deviations for random walk in random scenery.
\newblock Preprint (2005). To appear in {\it Annales Inst.\
H.\ Poincar\'e: Prob.\ Stat.}


\smallskip

\bibitem[G\"ar77]{Ga77}
{\sc J.~G\"artner},
\newblock On large deviations from the invariant measure,
\newblock {\it Th.~Prob.~Appl.} {\bf 22}, 24--39 (1977).

\smallskip

\bibitem[GH99]{GH99}
{\sc J.~G\"artner} and {\sc F.~den~Hollander,}
\newblock Correlation structure of intermittency in the
parabolic Anderson model,
\newblock {\em Probab.\ Theory Relat.\ Fields} {\bf 114}, 1--54 (1999).

\smallskip





\chwk{\bibitem[GH93]{GH93}
{\sc A.\ Greven} and {\sc F.\ den Hollander},
\newblock A variational characterization of the speed of a one-dimensional self-repellent random walk,
\newblock {\it Ann.\ Appl.\ Probab.} {\bf 3}, 1067-1099 (1993).}

\smallskip

\bibitem[HKM05]{HKM04}
{\sc R.~van der Hofstad, W.~K\"onig} and {\sc P.~M\"orters},
\newblock {The universality classes in the parabolic Anderson model}.
\newblock Preprint (2005). To appear in {\it Commun.\ Math.\ Phys.}

\smallskip


\bibitem[KS91]{KS91}
{\sc I.\ Karatzas} and {\sc S.E.\ Shreve},
\newblock {\it Brownian Motion and Stochastic Calculus},
2nd edition.
\newblock Springer, New York (1991).

\smallskip

\bibitem[KKS75]{KKS75}
{\sc H.\ Kesten, M.V.\ Kozlov} and {\sc F.\ Spitzer},
\newblock {A limit law for random walk in random environment.}
\newblock {\it Compositio Math.} {\bf 30}, 145--168 (1975).

\smallskip

\bibitem[Kni63]{Kn63}
{\sc F.B.\ Knight,}
\newblock { Random walks and a sojourn density process of
Brownian motion},
\newblock {\it Trans.\ Amer.\ Soc.} {\bf 109}, 56--86 (1963).

\smallskip

\bibitem[LeJ87]{LJ87}
{\sc Y.~Le Jan,}
\newblock { Temps local et superchamps,}
\newblock {\it S\'eminaire des Probabilit\'es XXI\ }, Lecture Notes in Math., {\bf 1247},
Springer, Berlin, 176--190 (1987).

\smallskip



\bibitem[Lut83]{Lut83}
{\sc J.~M.~Luttinger},
\newblock The asymptotic evaluation of a class of path integrals. {I}{I}.
\newblock {\it J. Math. Phys.}, {\bf 24:8}, 2070--2073 (1983).

\smallskip

\bibitem[MS87]{MS87}
{\sc P.~March} and {\sc A.-S.~Sznitman,}
\newblock {Some connections between excursion theory and
the discrete Schr\"odinger equation with random potentials}.
\newblock {\it Probab.\ Theory Relat.\ Fields} {\bf 109}, 11--53 (1987).

\smallskip

\bibitem[MR96]{MarRos96}
{\sc M.~B. Marcus} and {\sc J.~Rosen},
\newblock Gaussian chaos and sample path properties of additive functionals of
  symmetric {M}arkov processes.
\newblock {\em Ann. Probab.} {\bf 24:3}, 1130--1177 (1996).


\smallskip

\bibitem[McK80]{McK80}
{\sc A.J. McKane,}
\newblock Reformulation of $n \rightarrow 0$ models using anticommuting scalar
  fields.
\newblock {\em Physics Lett. A} {\bf 76}:22 (1980).


\smallskip


\bibitem[PS80]{PaSo80}
{\sc G.~Parisi} and {\sc N.~Sourlas,}
\newblock Self-avoiding walk and supersymmetry.
\newblock {\em J. Physique. Lettres.} {\bf 41}, L403--L406 (1980).

\smallskip

\bibitem[Ray63]{Ra63}
{\sc D.~Ray,}
\newblock Sojourn times of diffusion processes,
\newblock {\it Illinois J.\ Math.} {\bf 7}, 615-630 (1963).

\smallskip

\bibitem[RY91]{RY91}
{\sc D.\ Revuz} and {\sc M.\ Yor,}
\newblock {\it Continuous Martingales and Brownian Motion}.
\newblock Springer, Berlin (1991).

\smallskip

\bibitem[She85]{Shep85}
{\sc P.~Sheppard,}
\newblock On the Ray-Knight Markov property of local times.
\newblock {\it J.\ London Math.\ Soc.} {\bf 31}, 377--384 (1985).

\smallskip

\bibitem[Sym69]{Symanzik69}
{\sc K.~Symanzik,}
\newblock {E}uclidean quantum theory.
\newblock In R.~Jost, editor, {\em Local Quantum Theory}. Academic Press, New
  York, London (1969).

 \smallskip

\bibitem[Toth96]{Toth96}
{\sc B.~T\'oth},
\newblock Generalized Ray-Knight theory and limit theorems for
self-interacting random walks on $\Z^1$.
\newblock {\em Ann. Probab.} {\bf 24:3}, 1324--1367 (1996).



\end{thebibliography}
\end{document}